\documentclass[12pt]{amsart}
\usepackage{amssymb}
\usepackage{mathrsfs}
\usepackage{amsfonts}

\usepackage{amscd}

\theoremstyle{definition}
\theoremstyle{remark}
\numberwithin{equation}{section}

\begin{document}

\begin{center}
\bigskip {\large Szeg$\ddot{o}$ projection and matrix Hilbert transform in Hermitean Clifford analysis}

\smallskip

\medskip

\bigskip

Min Ku, Daoshun Wang \\

\smallskip
Department of Computer Science and Technology, Tsinghua
University, \\Beijing, 100084, P.R. China\\
E-mail: kumin0844{\rm @}163.com, kumin0844{\rm @}gmail.com

%\bigskip

%John Ryan

%\medskip

%Department of Mathematics

%University of Arkansas

%Fayetteville, AR 72701

\smallskip
\end{center}

\bigskip

\noindent \textbf{Abstract}\ The simultaneous null solutions of the
two complex Hermitean Dirac operators are focused on in Hermitean
Clifford analysis, where the matrix Hilbert transform was presented
and proved to satisfy the analogous properties of the Hilbert
transform in classical analysis and in orthogonal Clifford analysis.
Under this setting we will introduce the Szeg$\ddot{o}$ projection
operator for the Hardy space of Hermitean monogenic functions
defined on a bounded subdomain of even dimensional Euclidean space,
establish the Kerzman-Stein formula which closely connects the
Szeg$\ddot{o}$ projection operator with the Hardy projection
operator onto the Hardy space of Hermitean monogenic functions
defined on a bounded subdomain of even dimensional Euclidean space,
and get the Szeg$\ddot{o}$ projection operator in terms of the Hardy
projection operator and its adjoint. Further we will give the
algebraic and geometric characterizations for the matrix Hilbert
transform to be unitary in Hermitean Clifford analysis.

\bigskip

\noindent \textbf{Keywords:}\hspace{2mm} Hermitean Clifford
analysis, Szeg$\ddot{o}$ projection, Matrix Hilbert transformation,
Hardy space

\bigskip

\noindent \textbf{MSC(2000)}:\hspace{2mm} 30G35, 30C40, 31A10,
31A25, 31B10

\bigskip

\section{Introduction}

\medskip

\noindent The Hilbert transform in one dimensional space and its
properties were mainly developed by Titchmarsch and Hardy though it
is named after David Hilbert. This transform, which plays an
important role in engineering science such as signal analysis,
naturally appears when considering the boundary behavior of the
Cauchy transform. The crucial formula connecting the boundary value
of the Cauchy transform and Hilbert transform is the well-known
Plemelj-Sokhotzki formula. The classical multidimensional analogue
of the Hilbert transform is a tensorial one, studying the Riesz
transforms for each of the Cartesian variables
separately$\big{(}$see reference e.g. $[1]\big{)}$. As opposed to
these tensorial approaches, the orthogonal Clifford
analysis$\big{(}$seen in references e.g. $[2,3,4]\big{)}$
essentially provided a natural framework for generalizing a lot of
classical results from complex analysis and harmonic analysis in the
plane to the higher dimensional case. The central tool is the Cauchy
transform which leads to the Plemelj-Sokhotzki type formula when
taking the boundary values. Also the properties of the corresponding
singular operator was studied in virtue of the function theoretic
methods$\big{(}$see reference e.g. $[5]$ or elsewhere$\big{)}$.

\noindent In $[6]$, Kerzman and Stein proved the fundamental
property of the Cauchy transform $Cf$ of
$f\in\mathbf{L}_{2}\big{(}\Sigma\big{)}$ where $\Sigma$ is the
smooth boundary of a bounded open domain $D$ in the plane. They
stated that the operator $A=C-C^{*}\big{(}$where $C^{*}$ is the
adjoint operator of the operator $C \big{)}$ is a compact infinitely
smooth operator on $\mathbf{L}_{2}\big{(}\Sigma\big{)}$ and the
Szeg$\ddot{o}$ projection $S$ and the Hardy projection $C$ of
$\mathbf{L}_{2}\big{(}\Sigma\big{)}$ onto the Hardy space
$H^{2}\big{(}\Sigma\big{)}$ are related by Kerzman-Stein formula
$\big{(}$see references e.g. $[6,7,8,9]\big{)}$. Moreover, they
showed that the disc is the only plane region on which the Hilbert
transform $H$ on $\mathbf{L}_{2}\big{(}\Sigma\big{)}$ is unitary. In
$[10,11,12]$, Bernstein, Calderbank, Delanghe and their
collaborators generalized the Kerzman-Stein formula to the higher
dimensional case. Furthermore Delanghe$\big{(}$seen in $[12]\big{)}$
characterized the unitariness of the Hilbert transform under the
setting of orthogonal Clifford analysis. More related results on the
Szeg$\ddot{o}$ kernel and the Hilbert transform in orthogonal
Clifford analysis can be also found in references e.g. $[13,14-16]$.

\noindent More recently, offering yet a refinement of the orthogonal
case, Hermitean Clifford analysis in references e.g. $[17-23]$
emerged as a new and successful branch of Clifford analysis. It
focuses on the simultaneous null solutions of the two complex
Hermitean Dirac operators, which are invariant under the action of
the unitary group and were first studied in references e.g.
$[17-19]$. The Cauchy integral formula for Hermitean monogenic
functions defined in even dimensional Euclidean space taking values
in the complex Clifford algebra $\mathbb{C}_{2n}$ was constructed in
the framework of circulant $\big{(}2\times2\big{)}$ matrix
functions, and at the same time the intimate relationship with
holomorphic function theory of several complex variables
$\big{(}$see references e.g. $[24,25]\big{)}$ was established by
Brackx, De Schepper, Sommen and so on $\big{(}$see $[20]\big{)}$.
The Hermitean Cauchy transform, which gave rise to the Hardy
projection to be skew in Hermitean Clifford analysis, and the
related decomposition problems of continuous functions were
discussed in $[21,22]$. The new Hilbert-like matrix operator was
revealed by the non-tangential boundary limits of the Hermitean
Cauchy transform and the analogues of characteristic properties of
the matrix Hilbert transform in classical analysis and in orthogonal
Clifford analysis were given in $[23]$. Much recent progress can be
also seen elsewhere. Under this setting it is natural to think of
the orthogonal Szeg$\ddot{o}$ projection. However up to the present,
as for as we know, it has not been studied. In the underlying paper,
based on $[19-20,23,25,12,6,14]$, we will first define an inner
product on the space of square integral circulant
$\big{(}2\times2\big{)}$ matrix functions defined on the boundary of
a bounded subdomain in even dimensional Euclidean space, and
introduce the Szeg$\ddot{o}$ projection operator to be orthogonal
for the Hardy space of Hermitean monogenic functions defined on a
bounded subdomain of even dimensional Euclidean space. Then we will
establish the Kerzman-Stein formula which are closely related to the
Szeg$\ddot{o}$ projection operator and the Hardy projection operator
onto the Hardy space of Hermitean monogenic functions defined on a
bounded subdomain of even dimensional Euclidean space, and present
the Szeg$\ddot{o}$ projection operator in explicit terms of the
Hardy projection operator and its adjoint. Lastly we will give the
algebraic and geometric characterizations for the matrix Hilbert
transform to be unitary in Hermitean Clifford analysis.

\noindent The paper is organized as follows. In section $2$, we
recall some basic facts about Hermitean Clifford analysis which will
be needed in the sequel. In section $3$, we will introduce the
Szeg$\ddot{o}$ projection operator to be orthogonal for the Hardy
space of Hermitean monogenic functions defined on a bounded
subdomain of even dimensional Euclidean space, establish the
Kerzman-Stein formula which closely connects the Szeg$\ddot{o}$
projection operator with the Hardy projection operator onto the
Hardy space of Hermitean monogenic functions defined on a bounded
subdomain of even dimensional Euclidean space, and present the
Szeg$\ddot{o}$ projection operator in explicit terms of the Hardy
projection operator and its adjoint in Hermitean Clifford analysis.
In the last section we will give the algebraic and geometric
characterizations for the matrix Hilbert transform to be unitary in
Hermitean Clifford analysis.

\bigskip

\section{\noindent Preliminaries and notations}

\smallskip

\noindent In this section we  recall some basic facts about Clifford
algebra and Hermitean Clifford analysis which will be needed in the
sequel. More details can be also seen in the references e.g.
$[2,4,26,27,28-31]$ and $[17-23,25]$.

\smallskip

\noindent Let $\big{\{}e_{1},e_{2},\cdots,e_{m}\big{\}}$ be an
orthogonal basis of the Euclidean space $\mathbb{R}^{m}$, let
$\mathbb{R}^{m}$ be endowed with a non-degenerate quadratic form of
signature $\big{(}0,m\big{)}$ and let $\mathbb{R}_{0,m}$ be the
$2^{n}-$dimensional real Clifford algebra constructed over
$\mathbb{R}^{m}$ with basis
\begin{eqnarray*}
\Big{\{}e_{\mathcal{A}}:\mathcal{A}=\big{\{}h_{1},\cdots,h_{r}\big{\}}\in
\mathcal {P}\mathcal{N},1\leq h_{1}<h_{r}\leq m\Big{\}},
\end{eqnarray*}
where $\mathcal{N}$ stands for the set
$\big{\{}1,2,\cdots,m\big{\}}$ and $\mathcal{P}\mathcal{N}$ denotes
for the family of all order-preserving subsets of $\mathcal{N}$. We
denote $e_{\emptyset}$ as $e_{0}$ and $e_{\mathcal{A}}$ as
$e_{h_{1}\cdots h_{r}}$ for
$\mathcal{A}=\big{\{}h_{1},\cdots,h_{r}\big{\}}\in
\mathcal{P}\mathcal{N}$. The product in $\mathbb{R}_{0,m}$ is
defined by
\begin{displaymath}
\left\{\begin{array}{ll}
e_{\mathcal{A}}e_{\mathcal{B}}=(-1)^{N(\mathcal{A}\cap
\mathcal{B})}(-1)^{P(\mathcal{A},
\mathcal{B})}e_{\mathcal{A}\Delta\mathcal{B}},&\textrm{if
$\mathcal{A},\mathcal{B}\in\mathcal{P}\mathcal{N}$},\\
\lambda\mu=\sum\limits_{\mathcal{A},\mathcal{B}\in\mathcal{P}\mathcal{N}}
\lambda_{\mathcal{A}}\mu_{\mathcal{B}}e_{\mathcal{A}}e_{\mathcal{B}},&
\textrm{if
$\lambda=\sum\limits_{\mathcal{A}\in\mathcal{P}\mathcal{N}}\lambda_{\mathcal{A}}e_{\mathcal{A}},
\mu=\sum\limits_{\mathcal{B}\in\mathcal{P}\mathcal{N}}\mu_{\mathcal{B}}e_{\mathcal{B}}$},
\end{array}\right.
\end{displaymath}
where $N(\mathcal{A})$ is the cardinal number of the set
$\mathcal{A}$, and $P\big{(}\mathcal{A},\mathcal{B}\big{)} =
\sum\limits_{j\in \mathcal{B}}P(\mathcal{A},j)$,with
$P\big{(}\mathcal{A},j\big{)} = N\big{\{}i:i\in\mathcal
{A},i>j\big{\}}$. It follows $e_{0}$ is the identity element, now
written as $1$ and that in particular
\begin{displaymath}
\left\{\begin{array}{ll} e^{2}_{i}=-1,& \textrm{if $i=1,2,\cdots,m$},\\
e_{i}e_{j}+e_{j}e_{i}=0,&\textrm{if $1\leq i<j\leq m,$}\\
e_{h_{1}}e_{h_{2}}\cdots e_{h_{r}}=e_{h_{1}h_{2}\cdots
h_{r}},&\textrm{if $1\leq h_{1}<h_{2}<\cdots<h_{r}\leq m.$}
\end{array}\right.
\end{displaymath}
Thus the real Clifford algebra $\mathbb{R}_{0,m}$ is a real linear,
associative, but non-commutative algebra.

\smallskip

\noindent Any Clifford number $a$ in $\mathbb{R}_{0,m}$ may thus be
written as $a=\sum\limits_{N(\mathcal {A})=k}a_{\mathcal
{A}}e_{\mathcal {A}},a_{\mathcal {A}}\in\mathbb{R}$, or still as
$a=\sum\limits_{N(\mathcal {A})=k}[a]_{k}$, where
$[a]_{k}=\sum\limits_{N(\mathcal {A})=k}e_{\mathcal
{A}}[a]_{\mathcal {A}}$ is the so-called $k-$vector part of a
$\big{(}k=0,1,2,\cdots,m\big{)}$. The Euclidean space
$\mathbb{R}^{m}$ is embedded in $\mathbb{R}_{0,m}$ by identifying
$\big{(}x_{1},x_{2},\cdots, x_{m}\big{)}$ with the Clifford vector
$\underline{X}$ given by
$$\underline{X}=\sum\limits_{j=1}^{m}e_{j}x_{j}.$$
The conjugation in $\mathbb{R}_{0,m}$ is defined  as follows:
\begin{eqnarray*}
\bar{a}=\sum\limits_{\mathcal{A}}a_{\mathcal{A}}\bar{e}_{\mathcal{A}},
\bar{e}_{\mathcal{A}}=(-1)^{\frac{k(k+1)}{2}}e_{\mathcal{A}},N(\mathcal{A})=k,
a_{\mathcal{A}}\in \mathbb{R}.
\end{eqnarray*}
and hence
\begin{eqnarray*}
\overline{a b}=\bar{b}\bar{a}, \mbox{for arbitrary $a,b\in
\mathbb{R}_{0,m}$}.
\end{eqnarray*}

\noindent Note that the square of a vector $\underline{X}$ is scalar
valued and equals the norm squared up to a minus sign
$\underline{X}^{2}=-\big{<}\underline{X},\underline{X}\big{>}=-|\underline{X}|^{2}$.
The dual of the vector $\underline{X}$ is the vector valued first
order differential operator
$$\partial_{\underline{X}}=\sum\limits_{j=1}^{m}e_{j}\partial_{x_{j}}$$
is called Dirac operator. It is precisely this Dirac operator which
underlies the notion of monogenicity of a function, a notion which
is the higher dimensional counterpart or holomorphy in the complex
plane. As the Dirac operator factorizes the Laplacian,
$\Delta_{m}=-\partial^{2}_{\underline{X}}$, monogenicity can be
regarded as a refinement of harmonicity. We refer to this setting as
the orthogonal case, since the fundamental group leaving the Dirac
operator $\partial_{\underline{X}}$ invariant is the special
orthogonal group SO$\big{(}m;\mathbb{R}\big{)}$, which is doubly
covered by the Spin(m) group of the Clifford algebra
$\mathbb{R}_{0,m}$. For this reason, the Dirac operator is called a
rotation invariant
operator.\\
When allowing for complex constants and moreover taking the
dimension to be even, say $m=2n$, the same set of generators as
above, $\big{\{}e_{1},e_{2},\cdots,e_{2n}\big{\}}$, still satisfying
the above defining relation, may in fact also product the complex
Clifford algebra $\mathbb{C}_{2n}$. As $\mathbb{C}_{2n}$ is the
complexification of the real Clifford algebra $\mathbb{R}_{0,2n}$,
i.e. $\mathbb{C}_{2n}=\mathbb{R}_{0,2n}\oplus i\mathbb{R}_{0,2n}$,
any complex Clifford number $\lambda\in \mathbb{C}_{2n}$ may be
written as $\lambda=a+ib,a,b\in\mathbb{R}_{0,2n}$, leading to the
Hermitean conjugation
$\lambda^{\dag}=\big{(}a+ib\big{)}^{\dag}=\bar{a}-i\bar{b}$, where
the bar denotes the usual Clifford conjugation in
$\mathbb{R}_{0,2n}$, i.e. the main anti-involution for which
$\bar{e}_{j}=-e_{j},j=1,2,\cdots,2n$. This Hermitean conjugation
leads to a Hermitean inner product and its associated norm on
$\mathbb{C}_{2n}$ given by
$\big{(}\lambda,\mu\big{)}=\big{[}\lambda^{\dag}\mu\big{]}_{0}$ and
$\big{|}\lambda\big{|}=\sqrt{\big{[}\lambda^{\dag}\lambda\big{]}_{0}}=\Big{(}\sum\limits_{\mathcal
{A}}\big{|}\lambda_{\mathcal {A}}\big{|}^{2}\Big{)}^{\frac{1}{2}}$.
The above framework will be referred to as the Hermitean Clifford
analysis, as opposed to traditional orthogonal Clifford one.
Hermitean Clifford analysis then focuses on simultaneous null
solutions of two Hermitean Dirac operators
$\partial_{\underline{Z}}$ and
$\partial_{\underline{Z}^{\dagger}}$, introduced as follows.\\
 One of the ways for introducing Hermitean Clifford analysis is
by considering the complex Clifford algebra $\mathbb{C}_{2n}$ and a
so-called complex structure on it, i.e. an
SO$\big{(}2n,\mathbb{R}\big{)}$-element J for which
$J^{2}=-\textbf{1}$ $\big{(}$see e.g. $[17-20]\big{)}$. More
specifically, J is chosen to act upon the generators
$e_{1},e_{2},\cdots,e_{2n}$ of the Clifford algebra as
\begin{eqnarray*}
J[e_{j}]=-e_{n+j}~ \mbox{and}~ J[e_{n+j}]=e_{j},~ j=1,2,\cdots,n.
\end{eqnarray*}
Let us recall that the main objects of the Hermitean setting are
then conceptually obtained by considering the projection operators
$\frac{1}{2}\big{(}\textbf{1}\pm iJ\big{)}$ and letting them act on
the corresponding protagonists of the orthogonal framework. First of
all, the so-called Witt basis elements
$\Big{\{}f_{j},f^{\dagger}_{j}\big{|}j=1,2,\cdots,n\Big{\}}$ for the
complex Clifford algebra $\mathbb{C}_{2n}$ are obtained through the
action of $\frac{1}{2}(\textbf{1}\pm iJ)$ on the orthogonal basis
elements $e_{j}$:\begin{eqnarray*}
&&f_{j}=\frac{1}{2}\big{(}\textbf{1}+ i
J\big{)}[e_{j}]=\frac{1}{2}\big{(}e_{j}-ie_{n+j}\big{)}~,j=1,2,\cdots,n,\\
&&f^{\dagger}_{j}=-\frac{1}{2}\big{(}\textbf{1}- i
J\big{)}[e_{j}]=-\frac{1}{2}\big{(}e_{j}+ie_{n+j}\big{)}~,j=1,2,\cdots,n.
\end{eqnarray*}
These Witt basis elements satisfy the Grassmann identities
\begin{eqnarray*}
f_{j}f_{k}+f_{k}f_{j}=f^{\dag}_{j}f^{\dag}_{k}+f^{\dag}_{k}f^{\dag}_{j}=0,
~j,k=1,2,\cdots,n,
\end{eqnarray*}
and the duality identities
\begin{eqnarray*}
f_{j}f^{\dag}_{k}+f^{\dag}_{k}f_{j}=\delta_{jk}, ~j,k=1,2,\cdots,n.
\end{eqnarray*}
Next we identify a vector
$\underline{X}=\big{(}X_{1},X_{2},\cdots,X_{2n}\big{)}=\big{(}x_{1},x_{2},\cdots,x_{n},y_{1},\cdots,y_{n}\big{)}$
in $\mathbb{R}^{2n}$ with the Clifford vector
$\underline{X}=\sum\limits^{n}_{j=1}\big{(}e_{j}x_{j}+e_{n+j}y_{j}\big{)}$
and we denote by $\underline{X}|$ the action of the complex
structure J on $\underline{X}$, i.e.
\begin{eqnarray*}
\underline{X}|=J[\underline{X}]=\sum\limits^{n}_{j=1}\big{(}e_{j}y_{j}-e_{n+j}x_{j}\big{)}.
\end{eqnarray*}
Note that the vectors $\underline{X}$ and $\underline{X}|$ are
orthogonal w.r.t. the standard Euclidean scalar product, which
implies that the Clifford vectors $\underline{X}$ and
$\underline{X}|$ are both anti-commutative. The Hermitean Clifford
variables $\underline{Z}$ and $\underline{Z}^{\dagger}$ then arise
through the action of the projection operators on the standard
Clifford vector $\underline{X}$:
\begin{eqnarray*}
&&\underline{Z}=\frac{1}{2}\big{(}\textbf{1}+i
J\big{)}[\underline{X}]=\frac{1}{2}\big{(}\underline{X}+i\underline{X}|\big{)},\\
&&\underline{Z}^{\dagger}=-\frac{1}{2}\big{(}\textbf{1}-i
J\big{)}[\underline{X}]=-\frac{1}{2}\big{(}\underline{X}-i\underline{X}|\big{)}.
\end{eqnarray*}
They can be rewritten in terms of the Witt basis elements as
\begin{eqnarray*}
\underline{Z}=\sum\limits^{n}_{j=1}f_{j}z_{j},~\mbox{and}~
\underline{Z}^{\dagger}=(\underline{Z})^{\dagger}=\sum\limits^{n}_{j=1}f^{\dagger}_{j}z^{c}_{j},
\end{eqnarray*}
where $n$ complex variables $z_{j}=x_{j}+iy_{j}$ have been
introduced, with complex conjugates
$z^{c}_{j}=x_{j}-iy_{j},j=1,2,\cdots,n$. Finally, the Hermitean
Dirac operators $\partial_{\underline{Z}}$ and
$\partial_{\underline{Z}^{\dagger}}$ are derived out of the
orthogonal Dirac operator $\partial_{\underline{X}}$:
\begin{eqnarray*}
&&\partial_{\underline{Z}^{\dagger}}=\frac{1}{4}\big{(}\textbf{1}+i
J\big{)}[\partial_{\underline{X}}]=\frac{1}{4}\big{(}\partial_{\underline{X}}+i\partial_{\underline{X}|}\big{)},\\
&&\partial_{\underline{Z}}=-\frac{1}{4}\big{(}\textbf{1}-i
J\big{)}[\partial_{\underline{X}}]=-\frac{1}{4}\big{(}\partial_{\underline{X}}-i\partial_{\underline{X}|}\big{)},
\end{eqnarray*}
where we have introduced
\begin{eqnarray*}
\partial_{\underline{X}|}=J[\partial_{\underline{X}}]
=\sum\limits^{n}_{j=1}\big{(}e_{j}\partial_{y_{j}}-e_{n+j}\partial_{x_{j}}\big{)}.
\end{eqnarray*}
In terms of the Witt basis elements, the Hermitean Dirac operators
are expressed as
\begin{eqnarray*}
\partial_{\underline{Z}}=\sum\limits^{n}_{j=1}f^{\dagger}_{j}\partial_{z_{j}}~
\mbox{and}~
\partial_{\underline{Z}^{\dagger}}=(\partial_{\underline{Z}})^{\dagger}=
\sum\limits^{n}_{j=1}f_{j}\partial_{z^{c}_{j}},
\end{eqnarray*}
involving the classical Cauchy-Riemann operators
$\partial_{z_{j}}=\frac{1}{2}\big{(}\partial_{x_{j}}-i\partial_{y_{j}}\big{)}$
and their complex conjugates
$\partial_{z^{c}_{j}}=\frac{1}{2}\big{(}\partial_{x_{j}}+i\partial_{y_{j}}\big{)}$
in the complex $z_{j}-$planes, $j=1,2,\cdots,n$.

\smallskip

\noindent The Hermitean Dirac operators $\partial_{\underline{Z}}$
and $\partial_{\underline{Z}^{\dagger}}$ are invariant under the
action of a  realization, denoted $\widetilde{U}(n)$, of the unitary
group in terms of the Clifford algebras $\big{(}$see e.g.
$[17,19]\big{)}$. The group $\widetilde{U}(n)\subset
\mbox{Spin}(2n)$ is given by
\begin{eqnarray*}
\widetilde{U}(n)=\Big{\{}s\in \mbox{Spin} (2n)\big{|}~ \exists~
\theta\geq0:\tilde{s}I=e^{-i\theta}I\Big{\}}
\end{eqnarray*}
its definition involving the self-adjoint primitive idempotent
$I=I_{1}I_{2}\cdots I_{n}$, with
$I_{j}=f_{j}f^{\dagger}_{j}=\frac{1}{2}\big{(}1-ie_{j}e_{n+j}\big{)},~j=1,2,\cdots,n$.

\smallskip

\noindent Finally observe for further use that the Hermitean vector
variables and Dirac operators are isotropic, i.e.
\begin{eqnarray*}
(\underline{Z})^{2}=\big{(}\underline{Z}^{\dagger}\big{)}^{2}=0~
\mbox{and}~
\big{(}\partial_{\underline{Z}}\big{)}^{2}=\big{(}\partial_{\underline{Z}^{\dagger}}\big{)}^{2}=0.
\end{eqnarray*}
Whence the Laplacian
$\Delta_{2n}=-\partial_{\underline{X}}^{2}=-\partial_{\underline{X}|}^{2}$
allows the decomposition
\begin{eqnarray*}
\Delta_{2n}=4\big{(}\partial_{\underline{Z}}\partial_{\underline{Z}^{\dagger}}
+\partial_{\underline{Z}^{\dagger}}\partial_{\underline{Z}}\big{)}
\end{eqnarray*}
and one also has that
\begin{eqnarray*}
\underline{Z}~\underline{Z}^{\dagger}+\underline{Z}^{\dagger}\underline{Z}
=\big{|}\underline{Z}\big{|}^{2}=\big{|}\underline{Z}^{\dagger}\big{|}^{2}=\big{|}\underline{X}\big{|}^{2}
=\big{|}\underline{X}|\big{|}^{2}.
\end{eqnarray*}

\smallskip

\noindent For further use, we introduce the Hermitean oriented
surface elements $d\sigma_{\underline{Z}}$ and
$d\sigma_{\underline{Z}^{\dagger}}$ as follows
\begin{eqnarray*}&&\varepsilon(\underline{Z})
=\frac{2}{w_{2n}}\frac{\underline{Z}}{|\underline{Z}|^{2n}}~\mbox{and}~
\varepsilon^{\dagger}(\underline{Z})
=\frac{2}{w_{2n}}\frac{\underline{Z}^{\dagger}}{|\underline{Z}|^{2n}},\\
&&d\sigma_{\underline{Z}}=\sum\limits^{n}_{j=1}f^{\dagger}_{j}\widehat{dz_{j}}
~\mbox{and}~d\sigma_{\underline{Z}^{\dagger}}=\sum\limits^{n}_{j=1}f_{j}\widehat{dz^{c}_{j}}.
\end{eqnarray*}
Explicitly,\begin{eqnarray*} &&d\sigma_{\underline{Z}}
=-\frac{1}{4}(-1)^{\frac{n(n+1)}{2}}(2i)^{n}\big{(}d\sigma_{\underline{X}}
-id\sigma_{\underline{X}|}\big{)},\\
&&d\sigma_{\underline{Z}^{\dagger}}
=-\frac{1}{4}(-1)^{\frac{n(n+1)}{2}}(2i)^{n}\big{(}d\sigma_{\underline{X}}
+id\sigma_{\underline{X}|}\big{)},\\
&&\varepsilon=-\big{(}E+iE|\big{)},~\varepsilon^{\dagger}=\big{(}E-i
E|\big{)},
\end{eqnarray*}
where $d \sigma_{\underline{X}}$ denotes the vector valued oriented
surface element and $d \sigma_{\underline{X}|}=J[d
\sigma_{\underline{X}}]$. They are explicitly given by means of the
following differential forms of order $2n-1$
\begin{eqnarray*}
&&d \sigma_{\underline{X}}
=\sum\limits^{n}_{j=1}\big{(}e_{j}(-1)^{j-1}\widehat{dx_{j}}
+e_{n+j}(-1)^{n+j-1}\widehat{dy_{j}}\big{)},\\
&&d
\sigma_{\underline{X}|}=\sum\limits^{n}_{j=1}\big{(}e_{j}(-1)^{n+j-1}\widehat{dy_{j}}
+e_{n+j}(-1)^{j}\widehat{dx_{j}}\big{)},
\end{eqnarray*}
with
\begin{eqnarray*}
&&\widehat{dx_{j}}=dx_{1}\wedge\cdots\wedge dx_{j-1}\wedge
dx_{j+1}\wedge\cdots \wedge dx_{n}
\wedge dy_{1}\wedge\cdots\wedge dy_{n},\\
&&\widehat{dy_{j}}=dx_{1}\wedge\cdots\wedge dx_{n}\wedge
dy_{1}\wedge\cdots \wedge dy_{j-1}\wedge dy_{j+1}\wedge\cdots\wedge
dy_{n}.
\end{eqnarray*}

\noindent We denote the outward pointing unit normal vector at $X$
on $\partial\Omega$ by $\nu(\underline{X})$, $dS(\underline{X})$
stands for the classical element on $\partial\Omega$, then
\begin{eqnarray*} d \sigma_{\underline{X}}
=\nu(\underline{X})dS(\underline{X}),
~d\sigma_{\underline{X}|}=\nu|(\underline{X})dS(\underline{X}).
\end{eqnarray*}

\smallskip

\noindent In this context the functions taking
$\mathbb{C}_{2n}$-valued defined on an open region $\Omega$ of
$\mathbb{R}^{2n}$ will be considered. The continuity, continuously
differentiability, $\mathbb{L}_{p}(1< p<+\infty)$-integral and so on
of the function $f=\sum\limits_{\mathcal {A}}f_{\mathcal
{A}}e_{\mathcal {A}}: \Omega(\subset\mathbb{R}^{2n})\rightarrow
\mathbb{C}_{2n}$ where $f_{\mathcal
{A}}:\Omega(\subset\mathbb{R}^{2n})\rightarrow \mathbb{C}$, the
space of which are denoted respectively by
$\mathscr{C}\big{(}\Omega,\mathbb{C}_{2n}\big{)}$,
$\mathscr{C}^{1}\big{(}\Omega,\mathbb{C}_{2n}\big{)}$,
$\mathbf{L}_{p}\big{(}\Omega,\mathbb{C}_{2n}\big{)}$ and so on, are
ascribed to each component $f_{\mathcal {A}}$ which are respectively
continuous, continuously differential, $\mathbb{L}_{p}$-integrable
and so on. A function $f(\underline{X})$ defined and differentiable
in an open region $\Omega$ of $\mathbb{R}^{2n}$ with its boundary
$\partial\Omega$ and taking values in $\mathbb{C}_{2n}$ is called
(left) monogenic in $\Omega$ if
$\partial_{\underline{X}}f(\underline{X})=0$.

\noindent We introduce the particular circulant
$\big{(}2\times2\big{)}$ matrices
\begin{displaymath}\left.\begin{array}{ll}
\mathbf{D}_{(\underline{Z},\underline{Z}^{\dagger})}=\left(\begin{array}{ll}
\partial_{\underline{Z}}~ \partial_{\underline{Z}^{\dagger}}\\
\partial_{\underline{Z}^{\dagger}}~\partial_{\underline{Z}}
\end{array}\right),~
(\mathbf{D}_{(\underline{Z},\underline{Z}^{\dagger})})^{\dagger}=\left(\begin{array}{ll}
\partial_{\underline{Z}^{\dagger}}~ \partial_{\underline{Z}}\\
\partial_{\underline{Z}}~\partial_{\underline{Z}^{\dagger}}
\end{array}\right),
\mathcal{E}=\left(\begin{array}{ll}
\varepsilon \ \varepsilon^{\dagger}\\
\varepsilon^{\dagger} \  \varepsilon
\end{array}\right),
\mathbf{\underline{\delta}}=\left(\begin{array}{ll}
\delta \ 0\\
0 \ \delta
\end{array}\right),
\end{array}\right.
\end{displaymath}
then $\mathbf{D}_{(\underline{Z},\underline{Z}^{\dagger})}\mathcal
{E}=\mathbf{\underline{\delta}}(\underline{Z})$,
 i.e. $\mathcal {E}$ is the fundamental
solution of $\mathbf{D}_{(\underline{Z},\underline{Z}^{\dagger})}$
$\big{(}$see e.g. $[17,18,19,20]\big{)}$.\\
In the same setting of circulant $\big{(}2\times2\big{)}$ matrices,
we consider the functions $L_{1},L_{2},L\in
\mathscr{C}^{1}\big{(}\Omega,\mathbb{C}_{2n}\big{)}$ and the
corresponding circulant $\big{(}2\times2\big{)}$ matrix functions in
the following
\begin{displaymath}\left.\begin{array}{ll}
\mathcal {L}^{1}_{2}=\left(\begin{array}{ll}
L_{1}\ L_{2}\\
L_{2}\ L_{1}
\end{array}\right)\
\mbox{and}\ \mathcal {L}_{0}=\left(\begin{array}{ll}
L \ 0\\
0 \ L
\end{array}\right).
\end{array}\right.
\end{displaymath}
In the following context the operations of matrices such as addition
and multiplication, and the operations between the complex numbers
and the matrices respectively keep to the operation rules of the
usual numerical matrices and of multiplication between the complex
numbers and the usual numerical matrices.

\noindent \textbf{Definition $2.1$.}\ { Suppose that $\mathcal
{L}^{1}_{2}(\mathcal {L}_{0})\in \mathscr
{C}^{1}\big{(}\Omega,\mathbb{C}_{2n}\big{)}$ which means that each
entry of $\mathcal {L}^{1}_{2}(\mathcal {L}_{0})$ belongs to
$\mathscr {C}^{1}\big{(}\Omega,\mathbb{C}_{2n}\big{)}$. $\mathcal
{L}^{1}_{2}(\mathcal {L}_{0})$ is called as (left) H-monogenic if
and only if it satisfies the following system
 \begin{eqnarray*}
 \mathbf{D}_{(\underline{Z},\underline{Z}^{\dagger})}\mathcal
 {L}^{1}_{2}=\mathbf{0}\bigg{(}\mathbf{D}_{(\underline{Z},\underline{Z}^{\dagger})}\mathcal
 {L}_{0}=\mathbf{0}\bigg{)},
 \end{eqnarray*}}
where $\mathbf{0}$ denotes the $\big{(}2\times2\big{)}$ matrix with
zero entries. Similarly, it is obvious in the following that
$\mathcal {L}^{1}_{2}(\mathcal {L}_{0})\in\mathscr
{C}\big{(}\partial\Omega,\mathbb{C}_{2n}\big{)},\mathbb{H}^{\mu}\big{(}\partial\Omega,\mathbb{C}_{2n}\big{)}$,
$\mathbf{L}_{p}\big{(}\partial\Omega,\mathbb{C}_{2n}\big{)}\big{(}1<
p<+\infty\big{)}$ and so on which mean each entry of $\mathcal
{L}^{1}_{2}(\mathcal {L}_{0})$ belongs to $\mathscr
{C}\big{(}\partial\Omega,\mathbb{C}_{2n}\big{)},\mathbb{H}^{\mu}\big{(}\partial\Omega,\mathbb{C}_{2n}\big{)}$,
$\mathbf{L}_{p}\big{(}\partial\Omega,\mathbb{C}_{2n}\big{)}$ and so on.\\
In the following we introduce
\begin{eqnarray*}
&&\underline{V}=\frac{1}{2}\big{(}\underline{Y}+i\underline{Y}|\big{)},
\underline{V}^{\dagger}=-\frac{1}{2}\big{(}\underline{Y}-i\underline{Y}|\big{)},\\
&&d V_{(\underline{Z},\underline{Z}^{\dagger})}=\big{(}dz_{1}\wedge
dz^{c}_{1}\big{)}\wedge\big{(}dz_{2}\wedge
dz^{c}_{2}\big{)}\wedge\cdots\wedge\big{(}dz_{n}\wedge
dz^{c}_{n}\big{)},
\end{eqnarray*}
where $d V_{(\underline{Z},\underline{Z}^{\dagger})}$ denote the
Hermitean volume element.

\noindent For the functions $L_{i}\in
\mathbf{L}_{p}\big{(}\partial\Omega,\mathbb{C}_{2n}\big{)}\big{(}1<
p<+\infty,i=1,2\big{)}$, we define the orthogonal Cauchy type
integrals as follows
\begin{eqnarray*}
&&\mathcal {C}[L_{i}](\underline{Y})\triangleq\big{(}\mathcal
{C}_{\partial\Omega}L_{i}\big{)}(\underline{Y})=\int_{\partial\Omega}E\big{(}\underline{X}-\underline{Y}\big{)}d
\sigma_{\underline{X}}L_{i}(\underline{X}), \underline{Y}\notin
\partial\Omega,\\
&&\mathcal
{C}|[L_{i}](\underline{Y})\triangleq\big{(}\mathcal{C}|_{\partial\Omega}L_{i})(\underline{Y}\big{)}
=\int_{\partial\Omega}E|\big{(}\underline{X}-\underline{Y}\big{)}d
\sigma_{\underline{X}|}L_{i}(\underline{X}), \underline{Y}\notin
\partial\Omega,
\end{eqnarray*}
which are well-defined$\big{(}$see references e.g. $[4,15]\big{)}$,
where
$$E(\underline{X})=\frac{1}{w_{2n}}\frac{\overline{\underline{X}}}{|\underline{X}|^{2n}},
E|(\underline{X})=\frac{1}{w_{2n}}\frac{\overline{\underline{X}}|}{|\underline{X}|^{2n}},$$
and $d \sigma_{\underline{X}},d \sigma_{\underline{X}|}$ as above.
Then for $\underline{Y}\notin
\partial\Omega$, $$\partial_{\underline{Y}}\mathcal
{C}[L_{i}](\underline{Y})=0,\partial_{\underline{Y}|}\mathcal
{C}|[L_{i}](\underline{Y})=0\big{(}i=1,2\big{)}.$$

\smallskip

\noindent For the functions $\mathcal {L}^{1}_{2},\mathcal
{L}_{0}\in \mathbf{L}_{p}(\partial\Omega,\mathbb{C}_{2n})$, the
Hermitean Cauchy type integrals are defined by
\begin{eqnarray}
&&[\mathbf{C}\mathcal
{L}^{1}_{2}](\underline{Y})=\int_{\partial\Omega}\mathcal
{E}\big{(}\underline{Z}-\underline{V}\big{)}d
\Sigma_{(\underline{Z},\underline{Z}^{\dagger})}\mathcal
{L}^{1}_{2}(\underline{X}), \underline{Y}\notin
\partial\Omega,\\
&&[\mathbf{C}\mathcal {L}_{0}](\underline{Y})
=\int_{\partial\Omega}\mathcal
{E}\big{(}\underline{Z}-\underline{V}\big{)}d
\Sigma_{(\underline{Z},\underline{Z}^{\dagger})}\mathcal
{L}_{0}(\underline{X}), \underline{Y}\notin
\partial\Omega,
\end{eqnarray}
where
\begin{displaymath}\left.\begin{array}{ll}
d\Sigma_{(\underline{Z},\underline{Z}^{\dagger})}=\left(\begin{array}{ll}
d \sigma_{\underline{Z}} \ -d \sigma_{\underline{Z}^{\dagger}}\\
-d \sigma_{\underline{Z}^{\dagger}} \ d \sigma_{\underline{Z}}
\end{array}\right)\mbox{with $
d \sigma_{\underline{Z}}$ and  $d \sigma_{\underline{Z}^{\dagger}}$
as above}.
\end{array}\right.
\end{displaymath}

In the following we introduce the vector space
\begin{eqnarray*}
\mathscr{L}_{2}\big{(}\partial\Omega\big{)}=\left\{\mathcal
{L}^{1}_{2}=\left(
\begin{array}{ll}L_{1}\ L_{2}\\
L_{2} \ L_{1}\end{array}\right)\Big{|} L_{1}, L_{2}\in
\mathbf{L}_{2}\big{(}\partial\Omega,\mathbb{C}_{2n}\big{)}\right\},
\end{eqnarray*}
on which, inspired by the inner product
$\big{<}\cdot,\cdot\big{>}_{\mathbf{L}_{2}}$ on
$\mathbf{L}_{2}\big{(}\partial\Omega,\mathbb{C}_{2n}\big{)}$ given
by
\begin{eqnarray*}
%&&\big{<}L_{1},L_{2}\big{>}_{\mathbf{L}_{2}}
%=\int_{\partial\Omega}L^{\dagger}_{1}(\underline{X})L_{2}(\underline{X})d
%S_{\underline{X}},\\
%\tau
\big{<}L_{1},L_{2}\big{>}_{\mathbf{L}_{2}}
=\bigg{[}\int_{\partial\Omega}L^{\dagger}_{1}(\underline{X})L_{2}(\underline{X})d
S_{\underline{X}}\bigg{]}_{0},
\end{eqnarray*}
where $[\cdot]_{0}$ denotes the scale part of any $\cdot$ in
$\mathbb{C}_{2n}$. We introduce the following bilinear form
\begin{eqnarray*}
&&\big{<}\cdot,\cdot\big{>}_{\mathscr{L}_{2}} :
\mathscr{L}_{2}\big{(}\partial\Omega\big{)}\times\mathscr{L}_{2}\big{(}\partial\Omega\big{)}\rightarrow
\mathbb{C},\\
&&\left<\left(
\begin{array}{ll}L_{1}\ L_{2}\\
L_{2} \ L_{1}\end{array}\right),\left(
\begin{array}{ll}K_{1}\ K_{2}\\
K_{2} \
K_{1}\end{array}\right)\right>_{\mathscr{L}_{2}}\mapsto\left.
\begin{array}{ll}%\tau
\big{<}L_{1},K_{1}\big{>}_{\mathbf{L}_{2}}+%\tau
\big{<}L_{2},K_{2}\big{>}_{\mathbf{L}_{2}}.
\end{array}\right.\hspace{3cm}
\end{eqnarray*}
%where $tr(\cdot)$ denotes the trace operator of the
%$\big{(}2\times2\big{)}$ matrix $(\cdot)$ and
Then by directly calculating, for arbitrary $\mathcal
{L}^{1}_{2},\mathcal {K}^{1}_{2},\mathcal
{G}^{1}_{2}\in\mathbf{L}_{2}\big{(}\partial\Omega,\mathbb{C}_{2n}\big{)}$
where $\mathcal {G}^{1}_{2}$ is defined similarly to $\mathcal
{L}^{1}_{2}$ as above and arbitrary $\lambda\in\mathbb{C}$, we can
check
\begin{eqnarray*}
&&(i)~\big{<}\mathcal {L}^{1}_{2},\lambda\mathcal
{K}^{1}_{2}+\mathcal {G}^{1}_{2}\big{>}_{\mathscr{L}_{2}}
=\lambda\big{<}\mathcal {L}^{1}_{2},\mathcal
{K}^{1}_{2}\big{>}_{\mathscr{L}_{2}}
+\big{<}\mathcal {L}^{1}_{2},\mathcal {G}^{1}_{2}\big{>}_{\mathscr{L}_{2}},\\
&&(ii)~\bigg{(}\big{<}\mathcal {L}^{1}_{2},\mathcal
{K}^{1}_{2}\big{>}_{\mathscr{L}_{2}}\bigg{)}^{\dag}=\big{<}\mathcal
{K}^{1}_{2},\mathcal
{L}^{1}_{2}\big{>}_{\mathscr{L}_{2}},\\
&&(iii)~\big{<}\mathcal {L}^{1}_{2},\mathcal
{L}^{1}_{2}\big{>}_{\mathscr{L}_{2}}\geq0~\mbox{and}~\big{<}\mathcal
{L}^{1}_{2},\mathcal {L}^{1}_{2}\big{>}_{\mathscr{L}_{2}}=0
~\mbox{if and only if}~\mathcal {L}^{1}_{2}=\mathbf{0},
\end{eqnarray*}
where the operator $(\cdot)^{\dagger}$ as above and $\mathbf{0}$
denotes the $\big{(}2\times2\big{)}$ matrix with zero entries.
 Therefore
$\big{<}\cdot,\cdot\big{>}_{\mathscr{L}_{2}}$ is a inner product,
which derives the norm on
$\mathscr{L}_{2}\big{(}\partial\Omega\big{)}$ by
\begin{eqnarray*}\big{\|}\mathcal
{L}^{1}_{2}\big{\|}=\sqrt{\big{<}L_{1},L_{1}\big{>}_{\mathbf{L}_{2}}+\big{<}L_{2},L_{2}\big{>}_{\mathbf{L}_{2}}}.
\end{eqnarray*}
Hence
$\Big{(}\mathscr{L}_{2}\big{(}\partial\Omega\big{)},\big{\|}\cdot\big{\|}\Big{)}$
is the Hilbert space which is different from the space of
$\mathscr{L}_{2}\big{(}\partial\Omega\big{)}$ in references e.g.
$[23,25]$. Under this setting, we have the following Lemma without
proof which was also stated in $[23,25]$ in the sense of different
topology. For convenience without confusion and ambiguity,
$\Big{(}\mathscr{L}_{2}\big{(}\partial\Omega\big{)},\big{\|}\cdot\big{\|}\Big{)}$
still denotes by $\mathscr{L}_{2}\big{(}\partial\Omega\big{)}$ in
the following context.

\smallskip

\noindent \textbf{Lemma $2.1$.\ } {\it\ Suppose that $\Omega$ is an
open bounded subset of $\mathbb{R}^{2n}$ with smooth boundary
$\partial\Omega$. The functions $[\mathbf{C}{\mathcal
{L}}^{1}_{2}](\underline{X})$ and $[\mathbf{C}{\mathcal
{L}}_{0}](\underline{X})$ are defined similarly to
$[\mathbf{C}{\mathcal {L}}^{1}_{2}](\underline{Y})$ and
$[\mathbf{C}{\mathcal {L}}_{0}](\underline{Y})$ as above. If the
functions $\mathcal {L}^{1}_{2}(\underline{X}),\mathcal
{L}_{0}(\underline{X})\in\mathbf{L}_{p}\big{(}\partial\Omega,\mathbb{C}_{2n}\big{)}$,
$\big{(}1<p<+\infty\big{)}$,
then for arbitrary $\underline{T}\in
\partial\Omega$,
\begin{eqnarray*}
&&(i)~\mbox{for arbitrary $\underline{X}\in\mathbb{R}^{2n}\backslash
\partial\Omega$},
\mathbf{D}_{(\underline{Z},\underline{Z}^{\dagger})}\mathcal
{L}^{1}_{2}(\underline{X})=\mathbf{0},~\mathbf{D}_{(\underline{Z},\underline{Z}^{\dagger})}\mathcal
{L}_{0}(\underline{X})=\mathbf{0},\qquad\qquad\\
&&~~~~ \mbox{i.e. $\mathcal {L}^{1}_{2}(\underline{X}),\mathcal
{L}_{0}(\underline{X})$ are both
H-monogenic},\\
&&(ii)[\mathbf{C}{\mathcal
{L}}^{1}_{2}]^{\pm}(\underline{T})\triangleq\lim_{\Omega^{\pm}\ni
\underline{X}\rightarrow \underline{T}} [\mathbf{C}{\mathcal
{L}}^{1}_{2}](\underline{X})=(-1)^{\frac{n(n+1)}{2}}\frac{(2i)^{n}}{2}
\biggl(\pm{\mathcal {L}}^{1}_{2}(\underline{T}) +[\mathbf
{H}{\mathcal {L}}^{1}_{2}](\underline{T})\biggr),\\
&&~~~~~[\mathbf{C}{\mathcal
{L}}_{0}]^{\pm}(\underline{T})\triangleq\lim\limits_{\Omega^{\pm}\ni
\underline{X}\rightarrow \underline{T}} [\mathbf{C}{\mathcal
{L}}_{0}](\underline{X})=(-1)^{\frac{n(n+1)}{2}}(2i)^{n}
\frac{1}{2}\biggl(\pm{\mathcal {L}}_{0}(\underline{T})+[\mathbf
{H}{\mathcal {L}}_{0}](\underline{T})\biggr),\\
&&(iii)~[\mathbf{C}{\mathcal {L}}^{1}_{2}]^{\pm}(\underline{T}) \in
\mathbf{L}_{p}\big{(}\partial\Omega,\mathbb{C}_{2n}\big{)}\Big{(}[\mathbf{C}{\mathcal
{L}}_{0}]^{\pm}(\underline{T})\in
\mathbf{L}_{p}\big{(}\partial\Omega,\mathbb{C}_{2n}\big{)}\Big{)},
\end{eqnarray*}
where the limits of $(ii)$ mean the the non-tangential limits and it
is the same in this context,
\begin{displaymath}
\mathbf{H}=\frac{1}{2}\left(\begin{array}{ll}
\mathcal{H}+\mathcal{H}| \ -\mathcal{H}+\mathcal{H}|\\
-\mathcal{H}+\mathcal{H}| \ \mathcal{H}+\mathcal{H}|
\end{array}\right),
\end{displaymath}
and
\begin{eqnarray*}
&&[\mathcal
{H}f](\underline{T})=\mbox{p.v.}2\int_{\partial\Omega}E\big{(}\underline{Y}-\underline{T}\big{)}
d\sigma_{\underline{Y}}f(\underline{Y}),\underline{T}\in
\partial\Omega,\\
&&[\mathcal
{H}|f](\underline{T})=\mbox{p.v.}2\int_{\partial\Omega}E|\big{(}\underline{Y}-\underline{T}\big{)}
d\sigma_{\underline{Y}|}f(\underline{Y}),\underline{T}\in
\partial\Omega,
\end{eqnarray*}
which are both Cauchy principle value integrals in the sense of
$\mathbb{L}_{p}\big{(}1<p<+\infty\big{)}$. When the variables are
omitted without confusion and ambiguity, for convenience $[\mathcal
{H}f](\underline{T}),[\mathcal {H}|f](\underline{T})$ are for short
of $\mathcal {H}f,\mathcal {H}|f$ respectively and it is similar in
the following context.

\bigskip

\upshape {

\section{\noindent Szeg$\ddot{o}$ projection}

\smallskip

\noindent  In this section, we will introduce the Szeg$\ddot{o}$
projection operator for the Hardy space of Hermitean monogenic
functions defined on a bounded subdomain of even dimensional
Euclidean space, establish the Kerzman-Stein formula which is
closely related to the Szeg$\ddot{o}$ projection and the Hardy
projection for the Hardy space of Hermitean monogenic functions
defined on a bounded subdomain of even dimensional Euclidean space,
and get the Szeg$\ddot{o}$ projection operator in explicit terms of
the Hardy projection operator and its adjoint in the setting of
Hermitean Clifford analysis.

\noindent In what follows we will consider the Hardy space
\begin{eqnarray*}
\mathbb{H}^{2}\big{(}\Omega\big{)}=\Big{\{}\mathcal
{L}^{1}_{2}:\Omega\rightarrow\big{(}\mathbb{C}_{2n}\big{)}^{2\times2}\big{|}
\mathbf{D}_{(\underline{Z},\underline{Z}^{\dagger})}\mathcal
{L}^{1}_{2}=\mathbf{0}~\mbox{and}~L_{1}|_{\partial\Omega},L_{2}|_{\partial\Omega}\in
\mathbf{L}_{2}\big{(}\partial\Omega,\mathbb{C}_{2n}\big{)}\Big{\}}
\end{eqnarray*}
and $\mathbb{H}^{2}\big{(}\partial\Omega\big{)}$ denotes the
$\mathscr{L}_{2}\big{(}\partial\Omega\big{)}$-closure of the set of
boundary values of elements of $\mathbb{H}^{2}\big{(}\Omega\big{)}$.
Then associating Lemma $2.1$, the Hermitean Cauchy transform
$\mathbf{C}$ maps
$\mathbf{L}_{2}\big{(}\partial\Omega,\mathbb{C}_{2n}\big{)}$ onto
$\mathbb{H}^{2}\big{(}\partial\Omega\big{)}$ for arbitrary $\mathcal
{L}^{1}_{2}\in\mathbf{L}_{2}\big{(}\partial\Omega,\mathbb{C}_{2n}\big{)}$,
which is skew and so-called the Hardy projection.

\smallskip

\noindent By the same argument in $[23]$, associating the definition
of the above $\mathbb{C}$-valued inner product on
$\mathscr{L}_{2}\big{(}\partial\Omega\big{)}$, we have the following
Lemma which is only stated without proof.

\smallskip

\noindent \textbf{Lemma $3.1$\ }{\it\ Suppose that
$\mathscr{L}_{2}\big{(}\partial\Omega\big{)},\underline{\nu}$ and
$\mathbb{H}^{2}\big{(}\partial\Omega\big{)}$ as above. Then
\begin{eqnarray*}
&&(i)~ \mathbf{H}^{2}=\mathbf{I},\\
&&(ii)~ \mathbf{H}^{*}=\underline{\nu}\mathbf{H}\underline{\nu},\\
&&(iii)~\mbox{for arbitrary $\mathcal
{L}^{1}_{2}\in\mathscr{L}_{2}\big{(}\partial\Omega\big{)}$,}
\mathbf{H}\mathcal {L}^{1}_{2}=\mathcal {L}^{1}_{2}~\mbox{if and
only if}~ \mathcal
{L}^{1}_{2}\in\mathbb{H}^{2}\big{(}\partial\Omega\big{)}, \\
 &&(iv)~\mathcal
{L}^{1}_{2}\in\mathscr{L}_{2}\big{(}\partial\Omega\big{)}
=\mathbb{H}^{2}\big{(}\partial\Omega\big{)}
\oplus\underline{\nu}\mathbb{H}^{2}\big{(}\partial\Omega\big{)}\Big{(}\mbox{w.r.t.
}\big{<}\cdot,\cdot\big{>}_{\mathscr{L}_{2}}\Big{)},
\end{eqnarray*}
where
\begin{displaymath}\
\underline{\nu}=\frac{1}{2}\left(\begin{array}{ll}
\nu+\nu| \ -\nu+\nu|\\
-\nu+\nu| \ \nu+\nu|
\end{array}\right).
\end{displaymath}}

\smallskip

\noindent \textbf{Remark $3.1$}{\it\ The same results in Lemma $3.1$
were gotten in $[23]$ with respect to
$\big{(}\mathbb{C}_{2n}\big{)}^{2\times2}$-valued inner product
which is different from the above $\mathbb{C}$-valued inner product
on $\mathscr{L}_{2}\big{(}\partial\Omega\big{)}$.}

\noindent  The orthogonal projection operator $\mathbf{S}$ of
$\mathscr{L}_{2}\big{(}\partial\Omega\big{)}$ onto
$\mathbb{H}^{2}\big{(}\partial\Omega\big{)}$, which is so-called the
Szeg$\ddot{o}$ projection operator, may be Hermitean monogenically
extended to $\mathbb{H}^{2}\big{(}\Omega\big{)}$ by
\begin{eqnarray*}
\mathscr{S}\big{(}\mathcal
{L}^{1}_{2}(\underline{X})\big{)}=\int_{\partial\Omega}\mathcal
{\underline{S}}_{\underline{X}}(\underline{Y})\mathcal
{L}^{1}_{2}(\underline{Y})d S_{\underline{Y}},
\end{eqnarray*}
where $\mathcal {\underline{S}}_{\underline{X}}(\underline{Y})$ is
so-called the Szeg$\ddot{o}$ kernel.\\
That is, for arbitrary $\underline{X}\in\Omega$},
\begin{eqnarray*}
\mathscr{S}\big{(}\mathcal {L}^{1}_{2}(\underline{X})\big{)}
=\int_{\partial\Omega}\mathcal
{\underline{S}}_{\underline{X}}(\underline{Y})\mathcal
{L}^{1}_{2}(\underline{Y})d S_{\underline{Y}} =\mathcal
{L}^{1}_{2}(\underline{X}).
\end{eqnarray*}

\smallskip

\noindent \textbf{Remark $3.2$}{\it\ Particularly when $\Omega=B(1)$
the unit ball centered at $0$ of $\mathbb{R}^{2n}$,
$\partial\Omega=S^{2n-1}$ the unit sphere of $\mathbb{R}^{2n}$ and
$\nu(\underline{W})=\underline{W},\nu|(\underline{W})=\underline{W}|$
for arbitrary $\underline{W}\in S^{2n-1}$. Then
\begin{eqnarray*}
\mathscr{L}_{2}\big{(}S^{2n-1}\big{)}=\mathbb{H}^{2}\big{(}S^{2n-1}\big{)}
\oplus\underline{\nu}|_{S^{2n-1}}\mathbb{H}^{2}\big{(}S^{2n-1}\big{)},
\end{eqnarray*}
where\begin{displaymath}\
\underline{\nu}|_{S^{2n-1}}=\frac{1}{2}\left(\begin{array}{ll}
\underline{W}+\underline{W}| \ -\underline{W}+\underline{W}|\\
-\underline{W}+\underline{W}| \ \underline{W}+\underline{W}|
\end{array}\right).
\end{displaymath}

\smallskip

Given the boundary data $\mathcal { L}^{1}_{2}\in
\mathbf{L}_{2}\big{(}S^{2n-1},\mathbb{C}_{2n}\big{)}$, find the
function $\mathcal {K}^{1}_{2}$ such that
\begin{displaymath}\left\{\begin{array}{ll}
\underline{\Delta}\mathcal
{K}^{1}_{2}(\underline{X})=\mathbf{0}, \underline{X}\in B(1),\\
\mathcal{K}^{1}_{2}(\underline{X})
=\mathcal{L}^{1}_{2}(\underline{X}),\underline{X}\in S^{2n-1},
\end{array}\right.\quad\Leftrightarrow\left\{\begin{array}{ll}
\left.\begin{array}{ll}\Delta_{2n}
K_{1}(\underline{X})=0, \underline{X}\in B(1),\\
K_{1}(\underline{W})=L_{1}(\underline{W}),\underline{W}\in S^{2n-1},\end{array}\right.\\
\left.\begin{array}{ll}\Delta_{2n}
K_{2}(\underline{X})=0, \underline{X}\in B(1),\\
K_{2}(\underline{W})=L_{2}(\underline{W}),\underline{W}\in
S^{2n-1},\end{array}\right.
\end{array}\right.\end{displaymath}
where $\mathcal {L}^{1}_{2}$ as above,
$\left.\begin{array}{ll}\mathcal{K}^{1}_{2}=\left(\begin{array}{ll}
K_{1}~ K_{2}\\
K_{2}~ K_{1}
\end{array}\right)
\end{array}\right.$is defined similarly to $\mathcal
{L}^{1}_{2}$
and$$\left.\begin{array}{ll}\underline{\Delta}=\left(\begin{array}{ll}
\Delta_{2n}~ &0\\
0~ &\Delta_{2n}
\end{array}\right).
\end{array}\right.$$
In virtue of $(iv)$ in Lemma $3.1$, we have
\begin{eqnarray*}\mathcal {L}^{1}_{2}=\mathcal
{G}^{1}_{2}+\underline{\nu}\mathcal {H}^{1}_{2},\end{eqnarray*}
where $\mathcal {G}^{1}_{2},\mathcal
{H}^{1}_{2}\in\mathbb{H}^{2}\big{(}S^{2n-1}\big{)}$ are
defined similarly to $\mathcal {L}^{1}_{2}$. \\
  Then the above Dirichlet problem exists the unique solution.
Moreover the solution is formulated in the following form
\begin{eqnarray*}
\mathcal {K}^{1}_{2}(\underline{X})=\widetilde{\mathcal{G}}^{1}_{2}
+\mathcal {X}\widetilde{\mathcal{H}}^{1}_{2}, \underline{X}\in B(1),
\end{eqnarray*}
where
$\widetilde{\mathcal{G}}^{1}_{2},\widetilde{\mathcal{H}}^{1}_{2}\in\mathbb{H}^{2}\big{(}B(1)\big{)}$
are Hermitean monogenic extension of $\mathcal {G}^{1}_{2},\mathcal
{H}^{1}_{2}$ respectively $\big{(}$i.e.$\mathcal
{G}^{1}_{2},\mathcal {H}^{1}_{2}$ are the non-tangential boundary
value limits of
$\widetilde{\mathcal{G}}^{1}_{2},\widetilde{\mathcal{H}}^{1}_{2}\big{)}$
and$\left.\begin{array}{ll}\mathcal {X}=\left(\begin{array}{ll}
\underline{X}+\underline{X}| \ -\underline{X}+\underline{X}|\\
-\underline{X}+\underline{X}| \ \underline{X}+\underline{X}|
\end{array}\right).
\end{array}\right.$}
\smallskip

\noindent In what follows, we introduce the matrix Kerzman-Stein
operator on $\mathscr{L}_{2}\big{(}\partial\Omega\big{)}$ by
\begin{displaymath}\left.\begin{array}{ll}
\mathcal {\underline{A}}=\frac{1}{2}\left(\begin{array}{ll}
\mathcal {A}+\mathcal {A}| \ -\mathcal {A}+\mathcal {A}|\\
-\mathcal {A}+\mathcal {A}| \ \mathcal {A}+\mathcal {A}|
\end{array}\right),
\end{array}\right.
\end{displaymath}
where $\mathcal {A}=\mathcal {C}-\mathcal {C}^{*}$ and $\mathcal
{A}|=\mathcal {C}|-\mathcal {C}|^{*}$ are both well-defined, with
$\mathcal {C}^{*}$ and $\mathcal {C}|^{*}$ respectively denoting the
adjoint operators of $\mathcal {C}$ and $\mathcal {C}|$ on the
Hilbert space of
$\mathbf{L}_{2}\big{(}\partial\Omega,\mathbb{C}_{2n}\big{)}$ given
by
\begin{eqnarray*}
&&\mathcal {C}^{*}=\frac{1}{2}\big{(}1+\nu\mathcal {H}\nu\big{)}
:\mathbf{L}_{2}\big{(}\partial\Omega,\mathbb{C}_{2n}\big{)}\rightarrow H^{2}\big{(}\partial\Omega\big{)},\\
&&\mathcal {C}|^{*}=\frac{1}{2}\big{(}1+\nu|\mathcal
{H}|\nu|\big{)}:\mathbf{L}_{2}\big{(}\partial\Omega,\mathbb{C}_{2n}\big{)}\rightarrow
H^{2}\big{(}\partial\Omega\big{)},
\end{eqnarray*}
with
\begin{eqnarray*}
H^{2}\big{(}\partial\Omega\big{)}=\Big{\{}L_{1}:\Omega\rightarrow\mathbb{C}_{2n}\big{|}
\partial_{\underline{X}}L_{1}=0%\big{(}\mbox{or}~\partial_{\underline{X}|}
%f_{1}=0\big{)}
~\mbox{and}~L_{1}|_{\partial\Omega}\in
\mathbf{L}_{2}\big{(}\partial\Omega,\mathbb{C}_{2n}\big{)}\Big{\}}
\end{eqnarray*}
and $\nu,\nu|,\mathcal {H}$ as above and $1$ being the identity
operator. More detail can be seen in $[10,12,11]$.

\smallskip

\noindent Applying Lemma $3.1$, we directly get\\

\noindent \textbf{Lemma $3.2$}\ {\it\ Suppose that $\mathcal {A}$
and $\mathcal {A}|$ as above. Then
\begin{eqnarray}
\mathcal {\underline{A}}=\mathbf{C}-\mathbf{C}^{*}
=\mathbf{H}-\mathbf{H}^{*},~\mbox{i.e.}~ \mathcal
{\underline{A}}=\mathbf{H}-\underline{\nu}\mathbf{H}\underline{\nu},
\end{eqnarray}
where $\mathbf{H}^{*}$ as above and
$\mathbf{C}^{*}=\frac{1}{2}\big{(}\mathbf{I}+\mathbf{H}^{*}\big{)}$
mean the adjoint operators of $\mathbf{H}$ and $\mathbf{C}$ on
$\mathscr{L}_{2}\big{(}\partial\Omega\big{)}$.}

\noindent \textbf{Theorem $3.1$}\ {\it\ Suppose that $\mathbf{S}$
and $\mathbf{C}$ as above. Then
\begin{eqnarray}
\mathbf{S}\big{(}\mathbf{I}+\mathcal
{\underline{A}}\big{)}=\mathbf{C},
\end{eqnarray}
where $\mathbf{I}$ denote $\big{(}2\times2\big{)}$ identity matrix
operator.}

\noindent \textbf{Proof}\ {\it\ Since the operator $\mathbf{S}$ is
orthogonal projection operator on the Hilbert space
$\mathscr{L}_{2}\big{(}\partial\Omega\big{)}$, applying the property
of the orthogonal operator on the Hilbert space$\big{(}$see
reference e.g. $[32]\big{)}$, then $\mathbf{S}^{*}$ is well-defined,
where $\mathbf{S}^{*}$ means the adjoint operator of $\mathbf{S}$.
Moreover, $\mathbf{S}$ is the self-adjoint operator on
$\mathscr{L}_{2}\big{(}\partial\Omega\big{)}$, that is,
$\mathbf{S}=\mathbf{S}^{*}$.\\
Then as operators from $\mathscr{L}_{2}\big{(}\partial\Omega\big{)}$
to $\mathbb{H}^{2}\big{(}\partial\Omega\big{)}$,
\begin{eqnarray*}
\mathbf{S}\mathbf{C}=\mathbf{C}~\mbox{and}~\mathbf{C}\mathbf{S}=\mathbf{S}.
\end{eqnarray*}
Applying the property of the adjoint operator on the Hilbert space
of $\mathscr{L}_{2}\big{(}\partial\Omega\big{)}\big{(}$see reference
e.g. $[32]\big{)}$, $\big{(}\mathbf{S}\mathbf{C}\big{)}^{*}$ is
well-defined and
$\big{(}\mathbf{S}\mathbf{C}\big{)}^{*}=\mathbf{C}^{*}\mathbf{S}^{*}$,
where $\mathbf{C}^{*}$ means the adjoint operator on
$\mathscr{L}_{2}\big{(}\partial\Omega\big{)}$.  Taking the adjoint
operators with respect to $\big{<}.,.\big{>}_{\mathscr{L}_{2}}$, we
have
\begin{eqnarray*}
\mathbf{C}^{*}\mathbf{S}=\big{(}\mathbf{S}\mathbf{C}\big{)}^{*}=\mathbf{C}^{*}
~\mbox{and}~\mathbf{S}\mathbf{C}^{*}=\big{(}\mathbf{C}\mathbf{S}\big{)}^{*}=\mathbf{S}.
\end{eqnarray*}
 Hence\begin{eqnarray*}
\mathbf{S}=\mathbf{S}\mathbf{C}-\mathbf{S}\mathbf{C}^{*}=\mathbf{S}-\mathbf{C}.
\end{eqnarray*}
Therefore
\begin{eqnarray*}\mathbf{S}\big{(}\mathbf{I}+\underline{\mathcal
{A}}\big{)}=\mathbf{C},
\end{eqnarray*}
where $\mathbf{I}$ denotes
$\big{(}2\times2\big{)}$ identity matrix operator.\\
The proof of the result completes.}

\noindent \textbf{Remark $3.3$}{\it\ Theorem $3.1$ characterizes the
relation between Hermitean Hardy projection operator and
Szeg$\ddot{o}$ projection operator, which is the generalization of
classical Kerzman-Stein formula in the setting of Hermitean Clifford
analysis.}

We define the matrix operator as
follows\begin{displaymath}\left.\begin{array}{ll}
\mathbf{B}=\frac{1}{2}\left(\begin{array}{ll} \big{(}1+\mathcal
{A}\big{)}^{-1}+\big{(}1+\mathcal {A}|\big{)}^{-1} &
-\big{(}1+\mathcal {A}\big{)}^{-1}+\big{(}1+\mathcal {A}|\big{)}^{-1}\\
-\big{(}1+\mathcal {A}\big{)}^{-1}+\big{(}1+\mathcal
{A}|\big{)}^{-1}& \big{(}1+\mathcal
{A}\big{)}^{-1}+\big{(}1+\mathcal {A}|\big{)}^{-1}
\end{array}\right),
\end{array}\right.
\end{displaymath}
where $1$ denotes the identity operator on
$\mathbf{L}_{2}\big{(}\partial\Omega,\mathbb{C}_{2n}\big{)}$.\\
 Applying Lemma $4.5$ in $[10]$, the operator $1+\mathcal {A}$ and
$1+\mathcal {A}|$ are invertible on
$\mathbf{L}_{2}\big{(}\partial\Omega,\mathbb{C}_{2n}\big{)}\big{(}$also
see references $[11,32]$ or elsewhere$\big{)}$, the matrix operator
$\mathbf{B}$ is well defined on
$\mathscr{L}_{2}\big{(}\partial\Omega\big{)}$.

\smallskip

\noindent \textbf{Theorem $3.2$}\ {\it\ Suppose that $\mathbf{S}$
and $\mathbf{C}$ as above. Then the Szeg$\ddot{o}$ projection
operator is explicitly formulated by
\begin{eqnarray}
 \mathbf{S}=\mathbf{C}\big{(}\mathbf{I}+\mathcal
{\underline{A}}\big{)}^{-1}.
\end{eqnarray}
where $\mathbf{I}$ denote $\big{(}2\times2\big{)}$ identity matrix
operator.}

\noindent \textbf{Proof}\ {\it\ Since the matrix operators
$\mathbf{I}+\mathcal {\underline{A}}$ and $\mathbf{B}$ as above,
calculating directly, we get
\begin{eqnarray*}\mathbf{I}=
\left.\begin{array}{ll} \left(\begin{array}{ll}
2+\mathcal {A}+\mathcal {A}| & -\mathcal {A}+\mathcal {A}|\\
-\mathcal {A}+\mathcal {A}| & 2+\mathcal {A}+\mathcal {A}|
\end{array}\right)\mathbf{M},
\end{array}\right.
\end{eqnarray*}
where \begin{eqnarray*}\mathbf{M}=\left(\begin{array}{ll}
\big{(}1+\mathcal {A}\big{)}^{-1}+\big{(}1+\mathcal {A}|\big{)}^{-1}
\
-\big{(}1+\mathcal {A}\big{)}^{-1}+\big{(}1+\mathcal {A}|\big{)}^{-1}\\
-\big{(}1+\mathcal {A}\big{)}^{-1}+\big{(}1+\mathcal
{A}|\big{)}^{-1}\ \big{(}1+\mathcal
{A}\big{)}^{-1}+\big{(}1+\mathcal {A}|\big{)}^{-1}
\end{array}\right).
\end{eqnarray*}
 Then \begin{eqnarray*}\big{(}\mathbf{I}+\mathcal
{\underline{A}}\big{)}\mathbf{B}=\mathbf{I},\end{eqnarray*} i.e. the
matrix operator $\mathbf{I}+\mathcal {\underline{A}}$ is invertible
and its inverse is given by
\begin{eqnarray*}\big{(}\mathbf{I}+\mathcal
{\underline{A}}\big{)}^{-1}=\mathbf{B}.\end{eqnarray*} So it follows
that
\begin{eqnarray*}
 \mathbf{S}=\mathbf{C}\big{(}\mathbf{I}+\mathcal
{\underline{A}}\big{)}^{-1}.
\end{eqnarray*}}

\bigskip

\section{\noindent Characterization of matrix Hilbert transform}

\smallskip

\noindent In this section, we will give the algebraic and geometric
characterizations for the matrix Hilbert transform to be unitary in
Hermitean Clifford analysis, which is analogous to the corresponding
characterization of the Hilbert transform in classical analysis and
orthogonal Clifford analysis.

\smallskip

\noindent In the sequel we introduce the following functions
\begin{eqnarray*}
&&\alpha(\underline{X})=\frac{1}{2}\Big{(}1+i\nu(\underline{X})\Big{)},\quad
\beta(\underline{X})=\frac{1}{2}\Big{(}1-i\nu(\underline{X})\Big{)},
\underline{X}\in
\mathbb{R}^{2n},\\
&&\alpha|(\underline{X})=\frac{1}{2}\Big{(}1+i\nu|(\underline{X})\Big{)},\quad
\beta|(\underline{X})=\frac{1}{2}\Big{(}1-i\nu|(\underline{X})\Big{)},
\underline{X}\in \mathbb{R}^{2n}.
\end{eqnarray*}
 By directly calculating, it is easy to obtain the Lemma as
follows. \\
\noindent \textbf{Lemma $4.1$}\ {\it\ Suppose that
$\alpha(\underline{X}),\alpha|(\underline{X})$ and
$\beta(\underline{X}), \beta|(\underline{X})$ as above. Then
\begin{eqnarray}
&&(i)~\alpha^{2}(\underline{X})=\alpha(\underline{X}),
\beta^{2}(\underline{X})=\beta(\underline{X}),\\
&&(ii)~\alpha(\underline{X})\beta(\underline{X})=0,~\alpha|(\underline{X})\beta|(\underline{X})=0,\nonumber\\
&&(iii)~\alpha|(\underline{X})=\alpha|^{\dagger}(\underline{X}),
\beta|(\underline{X})=\beta|^{\dagger}(\underline{X}),\nonumber\\
&&(iv)~\alpha(\underline{X})+\beta(\underline{X})=1.
\end{eqnarray}

\noindent Related results can be also found in $[12,16,25]$ or
monographs on Fourier analysis elsewhere.

\noindent In what follows we introduce matrix functions
\begin{displaymath}\left.\begin{array}{ll}
\underline{\alpha}=\frac{1}{2}\left(\begin{array}{ll}
\alpha+\alpha| \ -\alpha+\alpha|\\
-\alpha+\alpha| \ \alpha+\alpha|
\end{array}\right),\
\underline{\beta}=\frac{1}{2}\left(\begin{array}{ll}
\beta+\beta| \ -\beta+\beta|\\
-\beta+\beta| \ \beta+\beta|
\end{array}\right).
\end{array}\right.
\end{displaymath}
where $\alpha,\beta,\alpha|,\beta|$ are for short of
$\alpha(\underline{X}),\beta(\underline{X}),\alpha|(\underline{X}),\beta|(\underline{X})$.
In the following context when without confusion and ambiguity, the
independent variable of considered functions are omitted.

\noindent  Making use of the above Lemma $4.1$ and directly
calculating of the matrix functions, we get the following Lemma.
\smallskip

\noindent \textbf{Lemma $4.2$} \ {\it\ Suppose that
$\underline{\alpha}$ and $\underline{\beta}$ as above. Then
\begin{eqnarray}
&&(i)~\underline{\alpha}\underline{\beta}
=\mathbf{0},\\
&&(ii)~\underline{\alpha}=\underline{\alpha}^{\dagger},
\underline{\beta}=\underline{\beta}^{\dagger},\nonumber\\
&&(iii)~\underline{\alpha}+\underline{\beta}=\mathbf{1},\underline{\nu}^{2}=-\mathbf{1},\\
&&(iv)\underline{\alpha}^{2}=\underline{\alpha},~\underline{\beta}^{2}=\underline{\beta},
\end{eqnarray}
where $\mathbf{0},\mathbf{1}$ denote $\big{(}2\times2\big{)}$ zero
matrix and identity matrix respectively.

\smallskip

\noindent Associating Lemmas $4.2,4.1$, we directly get the
algebraic decomposition of
$\mathscr{L}_{2}\big{(}\partial\Omega\big{)}$ as follows.

\noindent \textbf{Corollary $4.1$} \ {\it\ }For arbitrary $\mathcal
{L}^{1}_{2},\mathcal
{K}^{1}_{2}\in\mathscr{L}_{2}\big{(}\partial\Omega\big{)}$,
\begin{eqnarray*}
&&(i)~\big{<}\underline{\alpha}\mathcal
{L}^{1}_{2},\underline{\beta}\mathcal {K}^{1}_{2}\big{>}_{\mathscr{L}_{2}}=0,\\
&&(ii)~\mathscr{L}_{2}\big{(}\partial\Omega\big{)}=\underline{\alpha}\mathscr{L}_{2}\big{(}\partial\Omega\big{)}
\oplus\underline{\beta}\mathscr{L}_{2}\big{(}\partial\Omega\big{)}\Big{(}\mbox{w.r.t.
}\big{<}\cdot,\cdot\big{>}_{\mathscr{L}_{2}}\Big{)}.
\end{eqnarray*}
\noindent \textbf{Remark $4.1$}\ {\it\ The above Corollary $4.1$
gives the algebraic decomposition of
$\mathscr{L}_{2}\big{(}\partial\Omega\big{)}$. The analogous results
can be found in $[25]$, based on which the unique solution to the
classical Dirichlet problem on the unit ball in Hermitean Clifford
analysis is explicitly expressed.}

\smallskip

\noindent \textbf{Theorem $4.1$}\ {\it\ Suppose that $\Omega$ be a
bounded open domain of $\mathbb{R}^{2n}$ with smooth boundary
$\partial\Omega$. Let $\mathcal {L}^{1}_{2}$ and $\mathbf{0}$ be as
above. Then the following are equivalent
\begin{eqnarray*}
&&(i)~\underline{\alpha}\big{[}\mathbf {H}\underline{\alpha}\mathcal
{L}^{1}_{2}\big{]}=\mathbf{0},\underline{\beta}\big{[}\mathbf
{H}\underline{\beta}\mathcal {L}^{1}_{2}\big{]}=\mathbf{0}~
\mbox{for arbitrary $\mathcal {L}^{1}_{2}\in
\mathscr{L}_{2}\big{(}\partial\Omega\big{)}$}, \\
 &&(ii)~\mathbf
{H}\underline{\alpha}\mathcal {L}^{1}_{2}=\underline{\beta}\mathcal
{L}^{1}_{2}~\mbox{for arbitrary $\mathcal {L}^{1}_{2}\in
\mathbb{H}^{2}\big{(}\partial\Omega\big{)}$},\\
&&(iii)~\mathbf {H}\underline{\beta}\mathcal
{L}^{1}_{2}=\underline{\alpha}\mathcal {L}^{1}_{2}~\mbox{for
arbitrary $\mathcal {L}^{1}_{2}\in
\mathbb{H}^{2}\big{(}\partial\Omega\big{)}$},\\
&&(iv)~\mathbf {H}\underline{\nu}\mathcal
{L}^{1}_{2}=-\underline{\nu}\mathcal {L}^{1}_{2}~\mbox{for arbitrary
$\mathcal {L}^{1}_{2}\in
\mathbb{H}^{2}\big{(}\partial\Omega\big{)}$},\\
&&(v)~\mathbf {H}~\mbox{is unitary},~\mbox{i.e.}~ \mathbf {H}\mathbf
{H}^{*}=\mathbf {H}^{*}\mathbf {H}
=\mathbf{I}~\mbox{with $\mathbf{I}$ being identity matrix operator},\\
&&(vi)~\mathcal {\underline{A}}=\mathbf{0},\\
&&(vii)~\Omega~\mbox{is a ball},\\
&&(viii)~\mathcal
{\underline{S}}_{\underline{X}}(\underline{Y})=\mathbf
{C}_{\underline{X}}(\underline{Y}),\\ &&\mbox{with $\mathcal
{\underline{S}}_{\underline{X}}(\underline{Y}),\mathbf
{C}_{\underline{X}}(\underline{Y})$ denoting the Szeg$\ddot{o}$
kernel and the Cauchy kernel respectively,}\\ &&\mbox{i.e. the
Szeg$\ddot{o}$ kernel and the Cauchy kernel coincide.}
\end{eqnarray*}}

\smallskip

\noindent \textbf{Proof}\ {\it\ $``(i)\Rightarrow(ii)"$. For
arbitrary $\mathcal {L}^{1}_{2}\in
\mathbb{H}^{2}\big{(}\partial\Omega\big{)}$, by $(iii)$ in Lemma
$3.1$, we have
\begin{eqnarray}\label{k1}
\mathbf {H}\mathcal {L}^{1}_{2}=\mathcal
{L}^{1}_{2}.~\mbox{i.e.}~\underline{\beta}\mathbf {H}\mathcal
{L}^{1}_{2}=\underline{\beta}\mathcal {L}^{1}_{2}~\mbox{and}~
\underline{\alpha}\mathbf {H}\mathcal
{L}^{1}_{2}=\underline{\alpha}\mathcal {L}^{1}_{2}.
\end{eqnarray}
Associating the condition $(i)$ and $(iii)$ in Lemma $4.2$, we get
\begin{eqnarray}
&&\underline{\beta}\mathbf {H}\mathcal
{L}^{1}_{2}-\underline{\beta}\big{[}\mathbf
{H}\underline{\beta}\mathcal
{L}^{1}_{2}\big{]}=\underline{\beta}\mathcal {L}^{1}_{2},\nonumber\\
&&\underline{\beta}\big{[}\mathbf {H}\underline{\alpha}\mathcal
{L}^{1}_{2}\big{]}=\underline{\beta}\mathcal {L}^{1}_{2}.
\end{eqnarray}
Making use of the condition $\underline{\alpha}\big{[}\mathbf
{H}\underline{\alpha}\mathcal {L}^{1}_{2}\big{]}=0$, in term of
$(iii)$ in Lemma $4.2$, one gets
\begin{eqnarray}
\mathbf {H}\underline{\alpha}\mathcal
{L}^{1}_{2}=\underline{\alpha}\big{[}\mathbf
{H}\underline{\alpha}\mathcal
{L}^{1}_{2}\big{]}+\underline{\beta}\big{[}\mathbf
{H}\underline{\alpha}\mathcal
{L}^{1}_{2}\big{]}=\underline{\beta}\mathcal {L}^{1}_{2}.
\end{eqnarray}

\noindent $``(ii)\Rightarrow(i)"$. For arbitrary $\mathcal
{L}^{1}_{2}\in \mathscr{L}_{2}\big{(}\partial\Omega\big{)}$, using
$(iv)$ in Lemma $3.1$, we have
\begin{eqnarray}\mathcal {L}^{1}_{2}=\mathcal
{G}^{1}_{2}+\underline{\nu}\mathcal {H}^{1}_{2},\end{eqnarray} where
$\mathcal {G}^{1}_{2},\mathcal
{H}^{1}_{2}\in\mathbb{H}^{2}\big{(}\partial\Omega\big{)}$ are
defined similarly to $\mathcal {L}^{1}_{2}$. Therefore
\begin{eqnarray*}\mathbf {H}\underline{\alpha}\mathcal
{L}^{1}_{2}=\mathbf {H}\underline{\alpha}\mathcal
{G}^{1}_{2}+\mathbf {H}\underline{\alpha}\underline{\nu}\mathcal
{H}^{1}_{2}.\end{eqnarray*} The condition $(ii)$ acts on $\mathcal
{G}^{1}_{2}\in\mathbb{H}^{2}\big{(}\partial\Omega\big{)}$,
associating $(i)$ in Lemma $4.2$, we get
\begin{eqnarray*}\underline{\alpha}\mathbf {H}\underline{\alpha}\mathcal
{L}^{1}_{2}=\underline{\alpha}\mathbf {H}\underline{\alpha}\mathcal
{G}^{1}_{2}+\underline{\alpha}\mathbf
{H}\underline{\alpha}\underline{\nu}\mathcal
{H}^{1}_{2}=\underline{\alpha}\underline{\beta}\mathcal
{G}^{1}_{2}+\underline{\alpha}\mathbf
{H}\underline{\alpha}\underline{\nu}\mathcal
{H}^{1}_{2}=\underline{\alpha}\mathbf
{H}\underline{\alpha}\underline{\nu}\mathcal
{H}^{1}_{2}.\end{eqnarray*} i.e.
\begin{eqnarray*}i\underline{\alpha}\mathbf {H}\underline{\alpha}\mathcal
{L}^{1}_{2}=\underline{\alpha}\mathbf
{H}\underline{\alpha}i\underline{\nu}\mathcal
{H}^{1}_{2}.\end{eqnarray*} Applying the condition $(ii)$, we have
\begin{eqnarray*}\mathbf{0}=\underline{\alpha}\underline{\beta}\mathcal {H}^{1}_{2}
=\underline{\alpha}\mathbf {H}\underline{\alpha}\mathcal
{H}^{1}_{2}.\end{eqnarray*} Hence
\begin{eqnarray*}i\underline{\alpha}\mathbf {H}\underline{\alpha}\mathcal
{L}^{1}_{2}=\underline{\alpha}\mathbf
{H}\underline{\alpha}\big{(}\mathbf{1}+i\underline{\nu}\big{)}\mathcal
{H}^{1}_{2}=\underline{\alpha}\mathbf
{H}\underline{\alpha}^{2}\mathcal
{H}^{1}_{2}=\underline{\alpha}\mathbf {H}\underline{\alpha}\mathcal
{H}^{1}_{2}.\end{eqnarray*} By the condition $(ii)$, we get
\begin{eqnarray*}\underline{\alpha}\mathbf {H}\underline{\alpha}\mathcal
{L}^{1}_{2}=\underline{\alpha}\mathbf {H}\underline{\alpha}\mathcal
{H}^{1}_{2}=\underline{\alpha}\underline{\beta}\mathcal
{H}^{1}_{2}=0.\end{eqnarray*}

\noindent $``(ii)\Rightarrow(iii)"$. For arbitrary $\mathcal
{L}^{1}_{2}\in \mathbb{H}^{2}\big{(}\partial\Omega\big{)}$, applying
the term $(4.6)$ and the condition $(ii)$, we have
\begin{eqnarray*}
\mathbf {H}\underline{\beta}\mathcal {L}^{1}_{2}=\mathbf {H}\mathcal
{L}^{1}_{2}-\mathbf {H}\underline{\alpha}\mathcal
{L}^{1}_{2}=\mathcal {L}^{1}_{2}-\underline{\beta}\mathcal
{L}^{1}_{2}=\underline{\alpha}\mathcal {L}^{1}_{2}.
\end{eqnarray*}

\noindent $``(iii)\Rightarrow(ii)"$. It is similar to the procedure
of $(iii)\Rightarrow(ii)$.

\noindent $``(iii)\Rightarrow(iv)"$. For arbitrary $\mathcal
{L}^{1}_{2}\in \mathbb{H}^{2}\big{(}\partial\Omega\big{)}$, by the
condition $(iii)$ and the term $(4.6)$, we get
\begin{eqnarray*}
-i\mathbf {H}\underline{\nu}\mathcal {L}^{1}_{2}=2\mathbf
{H}\underline{\alpha}\mathcal {L}^{1}_{2}-\mathbf {H}\mathcal
{L}^{1}_{2}=2\underline{\alpha}\mathcal {L}^{1}_{2}-\mathcal
{L}^{1}_{2}=i\underline{\nu}\mathcal {L}^{1}_{2}.
\end{eqnarray*}
Hence\begin{eqnarray*} \mathbf {H}\underline{\nu}\mathcal
{L}^{1}_{2}=-\underline{\nu}\mathcal {L}^{1}_{2}.
\end{eqnarray*}

\noindent $``(iv)\Rightarrow(iii)"$. Since the procedure of
$(iii)\Rightarrow(iv)$, the result of $(iii)$ follows.

\noindent $``(iv)\Rightarrow(v)"$. For arbitrary $\mathcal
{L}^{1}_{2}\in \mathscr{L}_{2}\big{(}\partial\Omega\big{)}$, using
$(ii)$ in Lemma $3.1$ and the term $(4.9)$, we have
\begin{eqnarray*}\mathbf {H}\mathbf {H}^{*}\mathcal
{L}^{1}_{2}=\mathbf {H}\underline{\nu}\mathbf
{H}\underline{\nu}\mathcal {G}^{1}_{2}+\mathbf
{H}\underline{\nu}\mathbf {H}\underline{\nu}^{2}\mathcal
{H}^{1}_{2}.\end{eqnarray*} Making use of the condition $(iv)$, we
get \begin{eqnarray*}\mathbf {H}\mathbf {H}^{*}\mathcal
{L}^{1}_{2}=-\mathbf {H}\underline{\nu}^{2}\mathcal
{G}^{1}_{2}-\mathbf {H}\underline{\nu}\mathbf {H}\mathcal
{H}^{1}_{2}=\mathbf {H}\mathcal {G}^{1}_{2}-\mathbf
{H}\underline{\nu}\mathcal {H}^{1}_{2}=\mathcal
{G}^{1}_{2}+\underline{\nu}\mathcal {H}^{1}_{2}=\mathcal
{L}^{1}_{2}.\end{eqnarray*} Therefore for arbitrary $\mathcal
{L}^{1}_{2}\in \mathscr{L}_{2}\big{(}\partial\Omega\big{)}$,
\begin{eqnarray*}\mathbf {H}\mathbf {H}^{*}\mathcal
{L}^{1}_{2}=\mathcal {L}^{1}_{2}.~\mbox{i.e.}~\mathbf {H}\mathbf
{H}^{*}=\mathbf{I},\end{eqnarray*} where $\mathbf{I}$ denotes the
$\big{(}2\times2\big{)}$ identity matrix operator.

\noindent $``(v)\Rightarrow(iv)"$. Since $\mathbf {H}$ is unitary,
then $\mathbf{I}=\mathbf {H}^{*}\mathbf {H}.$ Associating $(ii)$ in
Lemma $3.1$, for arbitrary $\mathcal {L}^{1}_{2}\in
\mathbb{H}^{2}\big{(}\partial\Omega\big{)}$, we have
\begin{eqnarray*}\mathcal
{L}^{1}_{2}=\mathbf {H}^{*}\mathbf {H}\mathcal
{L}^{1}_{2}=\underline{\nu}\mathbf {H}\underline{\nu}\mathbf
{H}\mathcal {L}^{1}_{2}=\underline{\nu}\mathbf
{H}\underline{\nu}\mathcal {L}^{1}_{2}.\end{eqnarray*} Hence for
arbitrary $\mathcal {L}^{1}_{2}\in
\mathbb{H}^{2}\big{(}\partial\Omega\big{)}$,
\begin{eqnarray*}
\mathbf {H}\underline{\nu}\mathcal {L}^{1}_{2}=-
\underline{\nu}\mathcal {L}^{1}_{2}.\end{eqnarray*}

\noindent $``(v)\Rightarrow(vi)"$. From the condition $(v)$,
$\mathbf{I}=\mathbf {H}^{*}\mathbf {H}.$ Applying $(i)$ in Lemma
$3.1$, we have
\begin{eqnarray*}
\mathbf{H}=\mathbf {H}^{*}\mathbf {H}^{2}=\mathbf
{H}^{*}.\end{eqnarray*} By Lemma $3.2$, the result of $(vi)$
establishes.

\noindent $``(vi)\Rightarrow(v)"$. Applying Lemma $3.2$ and Lemma
$3.1$, it is easy to get the result.

\noindent $``(vi)\Rightarrow(vii)"$. From the condition $(vi)$,
$\mathcal {\underline{A}}=\mathbf{0}$. Hence $\mathcal
{\underline{A}}=\mathbf{H}-\underline{\nu}\mathbf{H}\underline{\nu}=\mathbf{0}$.
Since  \begin{displaymath}
\mathbf{H}=\frac{1}{2}\left(\begin{array}{ll}
\mathcal{H}+\mathcal{H}| \ -\mathcal{H}+\mathcal{H}|\\
-\mathcal{H}+\mathcal{H}| \ \mathcal{H}+\mathcal{H}|
\end{array}\right)~\mbox{and}~\underline{\nu}=\frac{1}{2}\left(\begin{array}{ll}
\nu+\nu| \ -\nu+\nu|\\
-\nu+\nu| \ \nu+\nu|
\end{array}\right),
\end{displaymath}
then \begin{displaymath}
\underline{\nu}\mathbf{H}\underline{\nu}=\frac{1}{8}\left(\begin{array}{ll}
\nu+\nu| \ -\nu+\nu|\\
-\nu+\nu| \ \nu+\nu|
\end{array}\right)\left(\begin{array}{ll}
\mathcal{H}+\mathcal{H}| \ -\mathcal{H}+\mathcal{H}|\\
-\mathcal{H}+\mathcal{H}| \ \mathcal{H}+\mathcal{H}|
\end{array}\right)\left(\begin{array}{ll}
\nu+\nu| \ -\nu+\nu|\\
-\nu+\nu| \ \nu+\nu|
\end{array}\right).
\end{displaymath}
Therefore
\begin{displaymath}\frac{1}{2}\left(\begin{array}{ll}
\mathcal{H}+\mathcal{H}| \ -\mathcal{H}+\mathcal{H}|\\
-\mathcal{H}+\mathcal{H}| \ \mathcal{H}+\mathcal{H}|
\end{array}\right)=\mathbf{H}=
\underline{\nu}\mathbf{H}\underline{\nu}=\frac{1}{2}\left(\begin{array}{ll}
\nu\mathcal{H}\nu+\nu|\mathcal{H}|\nu| \ -\nu\mathcal{H}\nu+\nu|\mathcal{H}|\nu|\\
-\nu\mathcal{H}\nu+\nu|\mathcal{H}|\nu| \
\nu\mathcal{H}\nu+\nu|\mathcal{H}|\nu|
\end{array}\right).
\end{displaymath}
Hence \begin{eqnarray*} &&\mathcal{H}=\nu\mathcal{H}\nu.
~\mbox{i.e.}~\nu_{Y}\big{(}Y-T\big{)}+\big{(}Y-T\big{)}\nu_{T}=0,\\
&&\big{<}Y-T,\nu_{T}\big{>}+\big{<}Y-T,\nu_{Y}\big{>}=0,
~\mbox{i.e.}~\big{<}Y-T,\nu_{T}+\nu_{Y}\big{>}=0,\end{eqnarray*}
where $Y,T\in\partial\Omega$ with $Y\neq T$ and $\nu_{Y},\nu_{T}$
denote outward pointing unit vectors at $Y,T\in\partial\Omega$
respectively $\big{(}$also see reference e.g. $[16]$ or
elsewhere$\big{)}$. The result $(vi)$ follows.

\noindent $``(vii)\Rightarrow(vi)"$. Since $\Omega$ is a ball and
the procedure of proof in ``$(vi)\Rightarrow(vii)$", the
$\mathcal{H}=\nu\mathcal{H}\nu$ and
$\mathcal{H}|=\nu|\mathcal{H}|\nu|$. So the result $(vi)$ holds.

\noindent  $``(vii)\Rightarrow(viii)"$. Since $\Omega$ is a
ball$\big{(}\mbox{i.e.}~\mathcal {\underline{A}}=\mathbf{0}\big{)}$,
by Theorem $3.2$, $\mathbf{S}=\mathbf{C}$. That is, the
Szeg$\ddot{o}$ projection and the Hardy projection coincide.

\noindent  $``(viii)\Rightarrow(vii)"$. If the Szeg$\ddot{o}$ kernel
and the Cauchy kernel coincide, $\mathbf{S}=\mathbf{C}$. As the
Szeg$\ddot{o}$ projection $\mathbf{S}$ is orthogonal on the Hilbert
space of $\mathscr{L}_{2}\big{(}\partial\Omega\big{)}$, then
$\mathbf{S}=\mathbf{S}^{*}$. Hence
$\mathbf{C}=\mathbf{S}=\mathbf{S}^{*}=\mathbf{C}^{*}$, where
$\mathbf{S}^{*},\mathbf{C}^{*}$ denote the adjoint operators of
$\mathbf{S},\mathbf{C}$ on
$\mathscr{L}_{2}\big{(}\partial\Omega\big{)}$. Then $\mathcal
{\underline{A}}=\mathbf{C}=\mathbf{C}^{*}=\mathbf{0}$ respectively.
That is, $(vii)$ establishes.

\noindent The proof of Theorem $4.1$ is complete.}

\smallskip

\noindent \textbf{Remark $4.2$} \ {\it\ ``$(vii)\Rightarrow(vi)$''
can be proved in virtue of direct calculation$\big{(}$see reference
$[25]\big{)}$, which leads the solutions to half Dirichlet problems
in the setting of Hermitean Clifford analysis.}

\noindent \textbf{Remark $4.3$}\ {\it\ The above theorem $4.1$
implies that the matrix Hilbert transform $\mathbf{H}$ is unitary if
and only if the bounded open subdomain $\Omega$ of $\mathbb{R}^{2n}$
is a ball. By Lemma $2.1$, the matrix Hilbert transform $\mathbf{H}$
is unitary if and only if the Hardy projection operator $\mathbf{C}$
is self-adjoint.}

\smallskip

\noindent {\bf{Acknowledgement}}\ {This research was supported in
part by Post-Doctor Foundation of China under Grant
$($No.$20090460316)$, the National Natural Science Foundation of
China under Grant $($No.$10871150$, No.$60873249)$, and $863$
Project of China under Grand $($No.$2008AA01Z419)$.}

}
\bigskip

\noindent {\bf References}

\smallskip

\upshape {

 \noindent[1] Stein E. M. Introduction to Fourier analysis on
Euclidean spaces. \emph{Princeton University Press}, Princeton, New
Jersey (1971)

\noindent [2] Brackx F., Delanghe R., Sommen F. Clifford analysis.
\emph{Res. Notes Math.}, {\bf 176}: Pitman, London (1982)
%\end{thebibliography}

\noindent{[3]} Delanghe R., Sommen F., Soucek V. Clifford algebra
and Spinor valued functions. \emph{Kluwer Academic}, Dordrecht
(1992)

\noindent{[4]} Gilbert J. E., Murry M. A. M. Clifford algebra and
Dirac operators in harmonic analysis. \emph{Cambridge Studies in
Advances Mathematics 26}, Cambridge University Press, Cambridge
(1991)

\noindent{[5]} McIntosh A. Clifford algebras, Fourier transforms,
and singular convolution operators on Lipschitz surfaces.
\emph{Revista Matem$\acute{ив}$tica Iberoamericana}, {\bf 10}(3):
665--721 (1994)

\noindent{[6]} Kerzman N., Stein E.M. The Cauchy kernel, the
Szeg$\ddot{o}$ kernel, and the Riemann mapping function \emph{Math.
Ann.}, {\bf 236}: 85--93 (1971)

\noindent{[7]} Bell S. The Cauchy transform, potential theory and
conformal mapping. \emph{CRC Press}, Boca Roton (1992)

\noindent{[8]} Bell S. Solving the Dirichlet problem in the plane by
means of the Cauchy integral. \emph{Indiana Univ. Math. Journal},
{\bf 39}(4): 1355--1371 (1990)

\noindent{[9]} Bell S. The Szeg$\ddot{o}$ projection and the
classical objects of potential theory in the plane. \emph{Duke
Mathematical Journal}, {\bf 64}(1): 1--26 (1991)

\noindent{[10]} Bernstein S., Lanzani L. Szeg$\ddot{o}$ projections
for Hardy spaces of monogenic functions and applications. \emph{Int.
J. Math. Math. Sc.}, {\bf 29}: 613--624 (2002)

\noindent{[11]} Calderbank D. Clifford analysis for Dirac operators
on manifolds with boundary. \emph{Max-Planck-Institute f$\ddot{u}$r
Mathematik}, Bonn (1996)

\noindent{[12]} Delanghe R. On some properties of the Hilbert
transform in Euclidean space. \emph{Bull. Belg. Math. Soc. Simon
Stevin}, {\bf 11}: 163--180 (2004)

\noindent{[13]} Constales D., Krausshar R.S. Szeg$\ddot{o}$ and
polymonogenic Bergman kernels for half-space and strip domains, and
single-periodic functions in Clifford analysis. \emph{Complex
Variables and Elliptic Equations}, {\bf 47}(4): 349--360 (2002)

\noindent{[14]} Delanghe R., Brackx F. Hypercomplex function theory
and Hilbert modules with reproducing kernel. \emph{Proc. London
Math. Soc.}, {\bf 37}(3): 545--576 (1978)

\noindent{[15]} Tao Qian. Hilbert transforms on the sphere and
Lipschitz surfaces. \emph{Quaternionic and Clifford Analysis},
Trends in Mathematics: 259--275 (2008)

\noindent{[16]} Delanghe R., Tao Q. Half Dirichlet problems and
decompositions of Poisson kernels. \emph{Advances in Applied
Clifford Algebras}, {\bf 17}(3): 383-393 (2007)

\noindent{[17]} Brackx F., Bures J., De Schepper H., Eelbode D.,
Sommen F., Souc$\check{e}$k V. Fundaments of Hermitean Clifford
analysis part I: complex structure. \emph{Complex Analysis and
Operator Theory}, {\bf 1}(3): 341--365 (2007)

\noindent{[18]} Brackx F., Bures J., De Schepper H., Eelbode D.,
Sommen F., Souc$\check{e}$k V. Fundaments of Hermitean Clifford
analysis part II: splitting of h-monogenic equations. \emph{Complex
Variables and Elliptic Equations}, {\bf 52} (10-11): 1063--1079
(2007)

\noindent{[19]} Brackx F., De Schepper H., Sommen F. The Hermitian
Clifford analysis toolbox. \emph{Advances in Applied Clifford
Algebras}, {\bf 18}(3-4): 451--487 (2008)

\noindent{[20]} Brackx F., De Knock B., De Schepper H., Sommen F. On
Cauchy and Martinelli-Bochner integral formulae in Hermitean
Clifford analysis. \emph{Bulletin of the Brazilian Mathematical
Society}, {\bf 40}(3): 395-416 (2009)

\noindent{[21]}  Abreu Blaya R., Bory Reyes J., Moreno Garc\'{i}a T.
Hermitian decomposition of continuous functions on a fractal
surface. \emph{Bull. Braz. Math. Soc.}, {\bf 40}(1): 107--115 (2009)

\noindent{[22]} Abreu Blaya R., Bory Reyes J., Brackx F., De Knock
Bram, De Schepper H., Pe\~{n}a Pe\~{n}a D., Sommen F. Hermitean
Cauchy integral decomposition of continuous functions on
hypersurfaces. \emph{Boundary Value Problems }, Article ID {\bf 57}:
425256, doi:10.1155 (2008)

\noindent{[23]} Brackx F., De Knock B., De Schepper H. A matrix
Hilbert transform in Hermitean Clifford analysis. \emph{J. Math.
Anal. Appl.}, {\bf 344}:1068--1078 (2008)

\noindent{[24]} Rocha-Chavez R., Shapiro M., Sommen, F. Integral
theorems for functions and differential forms in $\mathbb{C}_{m}$.
\emph{Research Notes in Mathematics 428,} Chapman Hall/CRC, New York
(2002)

\noindent{[25]} Min Ku, Daoshun Wang. Half Dirichlet problem for
matrix functions on the unit ball in Hermitean Clifford analysis.
\emph{J. Math. Anal. Appl.} (Submission for publication)

\noindent{[26]} Min Ku, Jinyuan Du. On integral representation of
spherical $k$-regular functions in Clifford analysis. \emph{Advances
in Applied Clifford Algebras}, {\bf 19}(1): 83--100 (2009)

\noindent{[27]} G$\ddot{u}$rlebeck K, Spr$\ddot{o}$ssig W.
Quaternionic and Clifford calculus for physicists and engineers.
\emph{Wiley: Chichester}, New York (1997)

\noindent{[28]} Min Ku. Integral formula of isotonic functions over
unbounded domain in Clifford analysis. \emph{Advances in Applied
Clifford Algebras}, {\bf 20}(1): 57--70 (2010)

\noindent{[29]} Ryan J. Complexified Clifford analysis.
\emph{Complex Var. Theory Appl.}, {\bf 1}(1): 119--149 (1982)

\noindent{[30]} Min Ku, Jinyuan Du, Daoshun Wang. On generalization
of Martinelli-Bochner integral formula using Clifford analysis.
\emph{Adv. Appl. Clifford Alge.}, {\bf 20}(2): 351--366 (2010)

\noindent{[31]} Min Ku, Jinyuan Du, Daoshun Wang. Some properties of
holomorphic Cliffordian functions in complex Clifford analysis.
\emph{Acta Mathematics Scientia}, {\bf 30B}(3): 747--768 (2010)

\noindent{[32]} Yosida K. Functional analysis (Fifth Edition).
\emph{Springer-Verlag} Berlin, Heidelberg, New York (1978)

}

\end{document}